\documentclass[10pt,a4paper]{article}
\usepackage{amsfonts}
\usepackage{amssymb}
\usepackage{amsmath}
\usepackage{epsfig}
\usepackage{amsthm}
\usepackage[latin1]{inputenc}
\usepackage[english]{babel}

\newtheorem{defi}{Definition}[section]
\newtheorem{proposition}[defi]{Proposition}
\newtheorem{theoreme}[defi]{Theorem}

\newtheorem{lemme}[defi]{Lemma}
\newtheorem{corollaire}[defi]{Corollary}
\newtheorem{remarque}[defi]{Remark}

\newcommand{\F}{\mathbb F} 
\newcommand{\Z}{\mathbb Z}
\newcommand{\N}{\mathbb N}

\newcommand{\p}{\mathbb P}

\newcommand{\id}{id}

\title{Large $p$-groups actions with a $p$-elementary abelian second ramification group.}
\author{Magali Rocher.}
\date{}

\begin{document}

\maketitle
 \begin{abstract}
Let $k$ be an algebraically closed field of characteristic $p>0$
and $C$ a connected nonsingular projective curve over $k$ with genus $g \geq 2$.
Let $(C,G)$ be a "big action" , i.e. a pair $(C,G)$ where $G$ is a $p$-subgroup of the $k$-automorphism group of $C$ such that$\frac{|G|}{g} >\frac{2\,p}{p-1}$. We denote by $G_2$ the second ramification group  of $G$ at the unique ramification point of the cover $C \rightarrow C/G$. The aim of this paper is to describe the big actions whose $G_2$ is $p$-elementary abelian. In particular, we obtain a structure theorem by considering the $k$-algebra generated by the additive polynomials.
We more specifically explore the case where there is a maximal number of jumps in the ramification filtration of $G_2$. In this case, we display some universal families.
\end{abstract}

\section{Introduction.}

\textit{Setting.}
Let $k$ be an algebraically closed field of characteristic $p>0$. We denote by $C$ a connected nonsingular projective curve over $k$, with genus $g \geq 2$, and by $G$ a $p$-subgroup of the $k$-automorphism group of $C$: $Aut_k(C)$, such that $\frac{|G|}{g}> \frac{2\,p}{p-1}$. Such a pair $(C,G)$ is called a "big action".   Then, there is a point of $C$ (say $\infty$) such that $G$ is equal to the wild inertia subgroup $G_1$ of $G$ at $\infty$ (cf. \cite{LM}). Moreover, the quotient curve $C/G$ is isomorphic to the projective line $\p^1_k$ and the ramification locus (respectively branch locus) of the cover $\pi: \; C \rightarrow C/G$ is the 
point  $\infty$ (respectively $\pi(\infty))$. 
Furthermore, the second lower ramification group $G_2$ of $G$ at $\infty$ is non trivial and is strictly included in $G_1$.
In addition, the quotient curve $C/G_2$ is isomorphic to $\p_k^1$ and the quotient group $G /G_2$ acts as a group of translations: $\{X \rightarrow X+y, \, y \in V\}$  of the affine line $C/G_2 - \{\infty \}$. 
\medskip

\noindent \textit{Motivation and purpose.}
When searching for a classification of big actions, it naturally occurs that the quotient $\frac{|G|}{g^2}$ has a "sieve" effect. Stichtenoth shows that, for any $p$-subgroup $G$ of $Aut_k(C)$, $\frac{|G|}{g^2} \leq \frac{4\,p}{(p-1)^2}$ (cf. \cite{St}). Later on, Lehr and Matignon prove that the big actions such that $\frac{|G|}{g^2} \geq \frac{4}{(p-1)^2}$ correspond to the $p$-cyclic \'etale covers of the affine line given by an Artin-Schreier equation: $W^p-W=f(X):=X\,S(X)+c\,X \in k[X]$, where $S(X)$ runs over the additive polynomials of $k[X]$ (cf. \cite{LM}). In a sequel paper \cite{MR3}, we go further in the classification and describe the big actions such that $\frac{|G|}{g^2} \geq \frac{4}{(p^2-1)^2}$. Under this condition, it is shown in \cite{MR} that $G_2$ is a $p$-elementary abelian group of order dividing $p^3$, hence the necessity to study big actions whose second ramification group $G_2$ is $p$-elementary abelian. 
\medskip

\noindent \textit{Outline of the paper.}
The main result of this paper is a structure theorem for the big actions with a $p$-elementary abelian second ramification group $G_2$. This result is obtained by considering the $k$-algebra generated by the additive polynomials of $k[X]$. Indeed, let $(C,G)$ be a big action whose $G_2$ is isomorphic to $(\Z/p\,\Z)^n$, with $n \geq 1$. Then, the function field of the curve is parametrized by $n$ Artin-Schreier equations: $W_i^p-W_i=f_i(X)\in k[X]$, with $1 \leq i \leq n$. For all $t\geq1$, we define $\Sigma_t$ as the $k$-subvector space of $k[X]$ generated by $1$ and the products of at most $t$ additive polynomials of $k[X]$. In section 3, we prove that each function $f_i$ belongs to $\Sigma_{i+1}$, which means that we can express $f_i$ as a linear combination over $k$ of products of at most $i+1$ additive polynomials of $k[X]$. This result generalizes the $p$-cyclic case, i.e. $n=1$, but, contrary to this case, the converse is no longer true, which means that such a family $(f_i)_{1\leq i\leq n}$ does not necessary give birth to a big action, except under specific conditions that are studied in what follows. The obstruction essentially lies in the embedding problem associated with the exact sequence:
$$ 0 \longrightarrow G_2 \simeq (\Z/p\, \Z)^n \longrightarrow G \longrightarrow V \longrightarrow 0$$
More precisely, we study the induced representation $\phi$: $G/G_2 \rightarrow Aut(G_2)\simeq Gl_n(\F_p)$ via the representation dual with respect to the Artin-Schreier pairing (see section 2). \\
\indent Sections 4 and 5 are devoted to two special cases of main interest. In section 4, we investigate the case where there is only one jump in the upper ramification filtration of $G_2$.  Then, the representation mentionned above is trivial or, equivalently, each function $f_i$ belongs to $\Sigma_2$. In section 5, we give a group-theoretic characterization of what can be regarded as the opposite case, namely: $f_i \in \Sigma_{i+1}-\Sigma_i$. Then, there is a maximal number of jumps in the upper ramificaton filtration of $G_2$. This case is relevant insofar as the representation $\phi$ is non trivial and provides much information.
To conclude, section 6 is devoted to examples illustrating section 5. In particular, we display a universal family parametrizing the big actions $(C,G)$ that satisfy $f_i \in \Sigma_{i+1}-\Sigma_i$, for $p=5$, a given $n \leq p-1$ and $dim_{\F_p} V=2$. When investigating the properties of the corresponding group $G$, we show that the center of $G$ is cyclic of order $p$ and relate $G$ with capable groups as defined by Hall (cf. \cite{Ha}).
\medskip

\noindent \textit{Notation and preliminary remarks.}
Let $k$ be an algebraically closed field of characteristic $p>0$.
We denote by $F$ the Frobenius endomorphism for a $k$-algebra. Then, $\wp$ means the Frobenius operator minus identity.
We denote by $k\{F\}$ the $k$-subspace of $k[X]$ generated by the polynomials $F^i(X)$, with $i \in \N$. It is a ring under the composition. Furthermore, for all $\alpha$ in $k$, $F\, \alpha=\alpha^p\,F$. The elements of $k\{F\}$ are the additive polynomials, i.e. the polynomials $P(X)$ of $k[X]$ such that for all $ \alpha$ and $\beta$ in $k$, $P(\alpha+ \beta) = P(\alpha)+ P(\beta)$. Moreover, a separable polynomial is additive if and only if the set of its roots is a subgroup of $k$ (see \cite{Go} chap. 1).\\
\indent Let $f(X)$ be a polynomial of $k[X]$. Then, there is a unique polynomial $red (f)(X)$ in $k[X]$, called the reduced representative of $f$, which is $p$-power free, i.e. $red(f)(X) \in \bigoplus_{(i,p)=1} k\, X^i$, and such that $red(f)(X)=f(X)$ mod $\wp (k[X]).$ We say that the polynomial $f$ is reduced mod $\wp(k[X])$ if and only if it coincides with its reduced representative $red(f)$.
The equation $W^p-W=f(X)$ defines a $p$-cyclic \'etale cover of the affine line that we denote by $C_f$. 
Conversely, any $p$-cyclic \'etale cover of the affine line $Spec \, k[X]$ corresponds to a curve $C_f$ where $f$ is a polynomial of $k[X]$ (see \cite{Mi} III.4.12, p. 127). 
By Artin-Schreier theory, the covers $C_f$ and $C_{red(f)}$ define the same $p$-cyclic covers of the affine line.
The curve $C_f$ is irreducible if and only if $red(f) \neq 0$.\\
\indent Throughout the text, the pair $(C,G)$ denotes a big action such that the second ramification group $ G_2$ is isomorphic to $(\Z / p \Z) ^n $, with $n \geq 1$.  We denote by $L:=k(C)$ the function field of $C$ and by  $k(X):=L^{G_2}$ the subfield of $L$ fixed by $G_2$. The extension $L/L^{G_2}$ is an \'etale cover of the affine line $Spec \, k[X]$ whose Galois group
  is $G_2 \simeq (\Z/p\, \Z)^n$. Therefore, it can be parametrized by $n$
  Artin-Schreier equations: $W^p-W=g_i(X)$, with $1 \leq i \leq n$. In other words,
$L=k(X,W_1,\cdots,W_n)$.
  As seen above, the functions $g_i(X)$ can be chosen in $k[X]$. Moreover, the quotient group
$G/G_2$ is a group of
  automorphisms of $k[X]$. Since it is a $p$-group, it actually acts as a group of translations of
Spec k[X], through $\tau_y: X  \rightarrow X+y$, where $y$ runs over a subgroup $V$ of $k$.  We remark that $V$ is an $\F_p$-subvector
  space of $k$. We denote by $v$ its dimension and thus obtain the exact sequence:
   $$ 0 \longrightarrow G_2 \simeq (\Z/p \,\Z)^n \longrightarrow G=G_1 
\stackrel{\pi}{\longrightarrow} V \simeq (\Z/\,p\, \Z)^v \longrightarrow 0$$ where for all $g$ in $G$,  $\pi(g):=g(X)-X$.
   We also fix a set theoritical section, i.e. a map $s: V \rightarrow G$, such that $\pi \circ s=
\id_V$.

\section{An embedding problem.}

\subsection{An $\F_p$-vector space dual of $G_2$.}
\indent We first exhibit an $\F_p$ - vector space dual of $G_2$.
 Following Artin-Schreier theory (see \cite{Bo}, chap. IX, ex. 19), we define the $\F_p$-vector space:
$$\tilde{A}:= \frac{ \wp (L) \cap k(X)}{ \wp (k(X))}$$
In other words, $\tilde{A}$ is the $\F_p$-vector space generated by the classes of the functions $g_i(X)$ modulo $\wp(k(X))$.
The inclusion $k[X] \subset k(X)$ induces an injection: 
$$ A:= \frac{ \wp (L) \cap k[X]}{ \wp (k[X])}\hookrightarrow  \tilde {A}$$ Since the extension is \'etale outside $\infty$, the functions $g_i(X)$ parametrizing the extension $L/k(X)$ can be chosen in $k[X]$. It follows that we can identify $A$ with $\tilde{A}$. Consider the Artin-Schreier pairing: 
$$ 
\left\{
\begin{array}{ll}
G_2 \times A  \longrightarrow \Z /p\, \Z \\
(g, \overline{\wp \, w}) \longrightarrow [g, \overline{\wp \, w} > := g(w)-w
\end{array}
\right.
$$
where $g$ belongs to $G_2 \subset Aut_k(L)$,  $w$ is an element of $L$ such that $\wp \, w \in k[X]$ and $\overline{\wp \, w}$ denotes the class of $\wp \, w$ mod $\wp (k[X])$. 
This pairing is non degenerate, which implies that,  as an $\F_p$-vector space, $A$ is dual to $G_2$.

\subsection{Two dual representations.}
\indent We now introduce two representations dual with respect to the Artin-Schreier pairing.
The first representation, say $\phi$, expresses the action of $G_1$ on $G_2$ via conjugation. 
Indeed, for all $y$ in $V$, we define an automorphism $\phi(y)$ of $G_2$ such that, for all $g$ in $G_2$, $\phi(y)(g):=s(y)^{-1} \, g\,  s(y)$.  Since $G_2$ is abelian, $\phi(y)$ does not depend on the lifting $s(y)$ in $G_1$ chosen for $y$ in $V$.
Therefore, there exists a representation $\phi$ which maps each $y$ in $V$ to $\phi(y)$ in $Aut(G_2)$.\\
\indent Then, we display a second representation expressing the action of $V$ on $A$.
More precisely, for all $y$ in $V$, we consider the automorphism $\rho(y)$ of $A$ defined as follows:
 $$  \rho(y): 
\left\{
\begin{array}{lc}
A \rightarrow A \\
\overline{\wp \, w} \rightarrow \overline{\wp (s(y)(w))}
\end{array}
\right.
$$ 
where $w$ is an element of $L$ such that $\wp \, w \in k[X]$.
As $s(y)$ belongs to $G_1 \subset Aut_k(L)$, then $s(y)(w)$ still lies in $L$. 
Furthermore, since $w^p-w \in k[X]$, then  $(s(y)(w))^p-s(y)(w)=s(y)(w ^p-w) \in s(y)(k[X]) \subset k[X]$, as $s(y)(X)=X+y$. This ensures that $\rho(y)$ is well-defined.
 Moreover, as $G_2$ trivially acts on $A \subset k(X)=L^{G_2}$, $\rho(y)$ is independent of the lifting
 $s(y)\in G_1$ chosen for $y$. Accordingly, we can define a representation $\rho$ which maps each $y$ in
 $V$ to $\rho(y)$ in $Aut(A)$.

\begin{remarque}
Note that for all $\overline{f(X)}$ in $A$ and for all $y$ in $V$, $\rho(y) \overline{f(X)}=\overline{f(X+y)}$.
\end{remarque}

\begin{proposition}
The two representations $\rho$ and $\phi$, as defined above, are dual with respect to the Artin-Schreier pairing.
\end{proposition}

\noindent \textbf{Proof:}
For all $y$ in $V$, for all $g$ in $G_2$ and for all $w$ in $L$ such that $\wp \, w$ is in $k[X]$,

$$
\begin{array}{rl}
 [ \phi(y)(g)\, , \,\overline{ \wp \, w} >& =
[s(y)^{-1} \; g \; s(y), \overline{\wp \, w} >\\
 &=s(y)^{-1}\, g\,s(y) (w)-w =s(y)^{-1}\, g \,s(y)(w)-s(y)^{-1}\; s(y) (w)\\
&=s(y)^{-1} \, (g \; s(y)(w)-s(y)(w)) = g \; s(y)(w)-s(y)(w)
\end{array}   
$$
since $g \; s(y)(w)-s(y)(w)= [g, \overline{\wp( s(y)(w))}> \in \F_p$.\\
As a conclusion, $[ \phi(y)(g) \,, \,\overline{ \wp \, w} >=[g \, , \, \overline{\wp (s(y)(w))} > =[g \,,\, \rho(y)(\overline{\wp \, w}) > \quad \square$
\medskip

Since the image of $\rho$ is a unipotent subgroup of $Gl_n(\F_p)$, one can find a basis for the $\F_p$-vector space $A$ in which the image of the representation $\rho$ can be identified with a subgroup of the upper triangular matrices in $Gl_n(\F_p)$. A means to do so is to endow $A$ with a filtration which proves to be dual of the upper ramification filtration of $G_2$.

\subsection{Dual filtrations on $A$ and $G_2$.}

 \indent The following three subsections are classical. Nevertheless, it is more convenient to recall both the proofs and the construction so as to fix the notation.
 
\subsubsection{A filtration and an adapted basis for $A$.}

\begin{defi}
\begin{enumerate}
\item We first gather from the canonical map "degree" a map defined on $A$ in the following way: 
$$  deg: 
\left\{
\begin{array}{ll}
A \rightarrow \N \cup \{ -\infty \} \\
 \overline{f(X)} \rightarrow \inf  \{ \, deg \; (f +\wp(P)), \, P \in k[X] \, \}
\end{array}
\right.
$$
\item For all $i$ in $\N$, we define a sequence of $\F_p$-subvector spaces of $A$ as follows: 
$$A^{i}:= \{ \overline{f(X)} \in A, \; deg (\overline{f(X)}) < i \}$$
\item From the increasing sequence: $\{0 \}=A^{0} \subset A^{ 1}\subset A^{ 2}\subset\cdots  A^{ r}\subset A^{
 r+1}=A$, we extract a strictly increasing sequence $(A^{\mu_i})_{0\leq i \leq s}$ such that:
$$\{0 \}=A^0=\cdots=A^{\mu_0}\subsetneq A^{ \mu_0+1}=\cdots=A^{\mu_1}\subsetneq A^{ \mu_1+1}= \cdots \subsetneq \cdots  A^{ \mu_{s}} \subsetneq A^{ \mu_s+1}=A$$
where the jumps $\mu_i$ are uniquely determined by the condition: $A^{\mu_i} \subsetneq A^{\mu_i+1}$. 
By definition of the function "degree" on $A$, all the integers $\mu_i$ are prime to $p$.
By convenience of notation, we also define $\mu_{s+1}$ as $\mu_s+1$ so that $A=A^{\mu_{s+1}}$.
For all $i$ in $\{0,\cdots,s+1\}$, we denote by $n_i$ the dimension of $A^{ \mu_i}$ over $\F_p$. Note that $n_0=0$ and $n_{s+1}=n$.
\item Starting from a basis of $A^{\mu_1}$, we complete it in a basis of $A^{ \mu_2}$, and so on until $A^{\mu_s+1}$. In this way, we construct a basis of $A$, say: $\{\overline{f_1(X)}, \cdots,\overline{f_n(X)} \}$, which is said to be "adapted" to the filtration defined above. Moreover, we impose specific conditions on the degree $m_i$ of each $\overline{f_i(X)}$:
\begin{enumerate}
\item  $\forall \; i \in \{1,\cdots,n\}, \; m_i$ is prime to $p$.
\item $\forall \; i \in \{1,\cdots,n-1\}, \; m_i \leq m_{i+1}$.
\item  $\forall \;  (\lambda_1,\cdots, \lambda_n) \in \F_p^n$ not all zeros, 
 $$ deg \; ( \sum_{i=1}^n \, \lambda_i  \, \overline{f_i(X)}) = \max_{i=1,\cdots,n} \{deg \;  \lambda_i \, \overline{f_i(X) } \}.$$
 \end{enumerate}
\end{enumerate}
 \end{defi}
 
\begin{remarque} Keeping the notation above, we notice that, for all $i \in \{0, \cdots, s\}$,
 $m_{n_i+1}=m_{n_i+2}=\cdots=m_{n_{i+1}}=\mu_{i}$.
\end{remarque}
 
This provides a new parametrization of the function field $L$. Indeed, for all $i$ in $\{1,\cdots,n\}$, we fix a representative mod $\wp(k[X])$ of $\overline{f_i(X)}$: $f_i(X)$ and  assume it to be reduced mod $\wp(k[X])$. As $m_i$ is prime to $p$, $f_i(X)$ still has degree $m_i$.
 We also suppose that for all $i$ in $\{1,\cdots,n\}$, $f_i(0)=0$. From now on, the extension $L/k(X)$ is parametrized by the $n$ Artin-Schreier equations: $W_i^p- W_i=f_i(X)$ with $1 \leq i \leq n$. 
 \medskip
 
 \subsubsection{The link with the upper ramification filtration of $G_2$.}

  In what follows, we highlight the correspondence between the jumps
 $(\mu_i)_{0 \leq i \leq s}$ in the filtration of $A$ and the jumps $(\nu_i)_{0\leq i \leq r}$ in the upper ramification filtration of $G_2$. Since $G_2$ is abelian, the Hasse-Arf Theorem (see e.g.  \cite{Se}, Chap. IV) asserts that the jumps in the upper ramification filtration are integers. So the ramification filtration reads as follows:
 $$G_2=(G_2)^0=\cdots =(G_2)^{\nu_0} \supsetneq (G_2)^{\nu_0+1} = \cdots =(G_2)^{\nu_1} \supsetneq \cdots
 =(G_2)^{\nu_r} \supsetneq (G_2)^{\nu_r+1}= \{0\}$$
 By convenience, put $\nu_{r+1}:=\nu_r+1$.

\begin{proposition}
 Keeping the notation above, $r=s$ and for all $i$ in $\{0,\cdots ,s+1\}$ , $\mu_i=\nu_i$. \\
 It follows that the filtration of $A$ and $G_2$ are dual with respect to the Artin-Schreier pairing, that is to say
 $(G_2)^{\nu_i}$ is the orthogonal of $A^{\mu_i}$, for all $i$ in $\{0,\cdots ,s+1\}$.
\end{proposition}
  
 \noindent \textbf{Proof:}
 Let $\nu_i$ be a jump in the upper ramification filtration of $G_2$, with $0 \leq i \leq r$. Since the $(G_2)^{\nu_i}$ are
 $\F_p$-subvectors spaces of $G_2$, one can find an index $p$-subgroup of $G_2$, say $\mathcal{H}$, such that
$(G_2)^{\nu_{i}+1} \subset \mathcal{H}$ and $(G_2)^{\nu_i} \not \subset \mathcal{H}$.
 As $L^{\mathcal{H}}/L^{G_2}$ is a $p$-cyclic cover of the affine line inside $L$, with Galois group equal to $G_2/\mathcal{H}$,
 it is parametrized by an Artin-Schreier equation: $W^p-W=f(X)=\sum_{i=1}^n \, \lambda_i \,f_i(X)$ with $(\lambda_i)_{1
 \leq i \leq n} \in (\F_p)^n-\{(0,0,\cdots,0)\}$.
Condition (c) in Definition 2.3.4 requires: 
 $ deg(f)= \max_{1 \leq i \leq n} \, \{deg \;  \lambda_i \, f_i(X)  \} \in \{m_i, \, 1 \leq i \leq n\}
 = \{\mu_i, \,0 \leq i \leq s\} $.
Besides, the group $G_2$ induces an upper ramification filtration on $G_2/\mathcal{H}$,
 namely $(\frac{G_2}{\mathcal{H}})^{\nu}=\frac{(G_2)^{\nu}\,\mathcal{H}}{\mathcal{H}}$ (see \cite{Se},  Chap. IV, Prop. 14).
 Therefore, the ramification filtration of $G_2/\mathcal{H}$ reads:
 $$  \Z/p\, \Z \simeq \frac{G_2}{\mathcal{H}}= (\frac{G_2}{\mathcal{H}})^0 = \cdots
 =(\frac{G_2}{\mathcal{H}})^{\nu_i}\supsetneq (\frac{G_2}{\mathcal{H}})^{\nu_i+1}=\{0\}$$
 This is precisely the $p$-cyclic case for which it is well-known 
 that the only jump of ramification: $\nu_i$ is equal to $deg(f)$ (see \cite{Se}, Chap. IV, ex. 4, p. 80). Therefore, $\nu_i
 \in \{\mu_j,\, 0 \leq j \leq s\} $. \\
 \indent Conversely, consider $\mu_i$, for $ 0 \leq i \leq s$.  Then, by Remark 2.4, $\mu_i=m_{n_{i+1}}$, i.e. the degree of the function $f_{n_{i+1}}$. The function field of the curve: $W^p-W=f_{n_{i+1}}(X)$, is a
 $p$-cyclic \'etale cover of the affine line whose Galois group is an index $p$-subgroup of $G_2$, say
 $\mathcal{H}$. We define the integer $\nu(G_2) \in \{\nu_i, \, 0 \leq i \leq r+1\}$ such that $G_2^{\nu(G_2)+1} \subset \mathcal{H}$ and $G_2^{\nu(G_2)} \not \subset \mathcal{H}$. As seen above, $m_{n_{i+1}}= \nu(G_2)$.
Therefore, $\mu_i \in \{\nu_j, \,0 \leq j \leq r\}$. Accordingly, $\{\nu_i, \, 0 \leq i \leq r\} = \{\mu_i,\, 0 \leq i \leq s\}$. They are both strictly increasing sequence, so $r=s$ and for all $i$ in $\{0,\cdots ,s\}$ , $\mu_i=\nu_i$. 
 In addition, $\mu_{s+1}=\mu_s+1=\nu_r+1=\nu_{r+1}$, which completes the proof of the proposition. $\square$

\subsubsection{The different exponent and the genus of the extension.}
In this section, we establish a formula to calculate the different exponent and the genus of the extension $L/L^{G_2}$. We keep the notation defined in sections 2.3.1 and 2.3.2.

\begin{proposition}
Let $(C,G)$ be a big action such that $G_2 \simeq (\Z/p\Z)^n$, with $n \geq 1$. 
The different exponent of the extension $L/L^{G_2}$ is given by the formula:
$$d= (p-1)\, \sum_{i=1}^n \, p^{i-1} \, (m_i+1)$$
\end{proposition}

\noindent \textbf{Proof:} 
Since $G_2$ is abelian, one can apply to $L/L^{G_2}$ the upper index version of the Hilbert's different formula as given in \cite{Au} (p. 120): $d= \sum_{i=0}^{\infty} \, (|G_2| - [G_2: (G_2)^i])$.
In our case, this formula reads:
$$d=  (\nu_0+1)\,(|G_2|-[G_2: (G_2)^{\nu_0}]) + \sum_{j=1}^{r}\, (\nu_{j}-\nu_{j-1})\, (|G_2|-[G_2:(G_2)^{\nu_j}])$$
Using Proposition 2.5, we obtain:
$$
\begin{array}{ll}
d&= (\nu_0+1)\,(|G_2|-|A^{\mu_0}|)+ \sum_{j=1}^{r} \, (\nu_{j}-\nu_{j-1})\,(|G_2|-|A^{\mu_j}|)\\
& \\
&=(\mu_0+1)\,  (p^n-p^{n_0})+\sum_{j=1}^{s} \, (\mu_{j}-\mu_{j-1})\, (p^n-p^{n_{j}})\\
&\\
&= \sum_{j=0}^{s} (p^n-p^{n_{j}})\, (\mu_{j}+1)+\sum_{j=1}^{s+1} (p^{n_{j}}-p^n)\, (\mu_{j-1}+1)\\
&\\
&= \sum_{j=1}^{s+1} (p^n-p^{n_{j-1}})\, (\mu_{j-1}+1)+\sum_{j=1}^{s+1} (p^{n_{j}}-p^n)\, (\mu_{j-1}+1)\\
&\\
&=  \sum_{j=1}^{s+1}  (p^{n_j}-p^{n_{j-1}})\, (\mu_{j-1}+1)\\
&\\
&= \sum_{i=1}^{s+1} \sum_{j=n_{i-1}+1}^{n_i} p^{j-1} \,(p-1) \, (\mu_{i-1}+1)\\
&\\
&=(p-1) \sum_{i=1}^{s+1} \sum_{j=n_{i-1}+1}^{n_i} p^{j-1} \, (m_j+1)\\
&\\
&= (p-1)\, \sum_{i=1}^n \, p^{i-1} \, (m_i+1) \qquad \quad \square
\end{array}
$$

Note that another proof of this formula can be obtained by applying the formula given by Garcia and Stichtenoth in \cite{GS}.

\begin{corollaire}
Let $(C,G)$ be a big action such that $G_2 \simeq (\Z/p\Z)^n$, with $n \geq 1$. 
The genus of the extension $L/L^{G_2}$ is given by the formula:
$$g= \frac{1} {2} \,(p-1)\, \sum_{i=1}^n \, p^{i-1} \, (m_i-1) $$
\end{corollaire}

\textbf{Proof}:
The formula directly derives from the Hurwitz genus formula (see e.g. \cite{St93}) combined with the formula given in Proposition 2.6:
$$2\, (g-1)= 2\, (g_{C/G_2} -|G_2|)+d=-2\,p^n+(p-1)\, \sum_{i=1}^n \, p^{i-1} \, (m_i+1)
\qquad \square$$

 \subsection{Matricial representations of $\rho$ and $\phi$.}
 \indent From now on, we work in the adapted basis constructed for $A$ in section 2.3.1: $\{\overline{f_1(X)}, \cdots, \overline{f_n(X)}\}$.
 For any $y$ in $V$, we denote by $L(y)$ the matrix of the automorphism $\rho(y)$ in this basis.
As indicated in Remark 2.1, we recall that for all $y$ in $V$ and for all $i$ in $\{1,\cdots,n\}$, $ \rho(y)\,  \overline{f_i(X)}= \overline{f_i(X+y)}$. Moreover, the conditions imposed on the degree of the functions $\overline{f_i(X)}$ imply that the matrix
 $L(y)$ belongs to $T_{1,n}^u(\F_p)$, the subgroup of $Gl_n(\F_p)$ made of the upper triangular matrices with identity
 on the diagonal. Thus, $L(y)$ reads as follows:
  $$ L(y):=
\begin{pmatrix}
 1 &\ell_{1,2}(y)&  \ell_{1,3}(y) &\cdots& \ell_{1,n}(y) \\
 0 & 1 &\ell_{2,3}(y) & \cdots&\ell_{2,n}(y) \\
 0 & 0 &\cdots&\cdots&\ell_{i,n}(y) \\
 0 & 0 & 0 & 1& \ell_{n-1,n}(y)\\
 0 & 0 & 0 & 0 & 1
\end{pmatrix} \in Gl_n(\F_p)
$$
In other words, 
$$ \forall \, y \in V, \;, f_1(X+y)-f_1(X)=0 \quad \mod \; \wp(k[X])$$
\begin{equation} \label{1}
\forall \, i \in \{2,\cdots,n\}, \; \forall \, y \in V, \;, f_i(X+y)-f_i(X)=\sum_{j=1}^{i-1}\,
\ell_{j,i}(y) \, f_j(X) \quad \mod \; \wp(k[X]) 
\end{equation}

\begin{remarque}
We still denote by $m_i$ the degree of the function $f_i$. We observe that the degree of the left-hand side of $\eqref{1}$ is at most $m_i-1$. It follows that, whenever $m_i=m_j$, $f_j$ does not occur in the right-hand side of $\eqref{1}$, which means that $\ell_{j,i}$ is zero on $V$. 
\end{remarque}

 \begin{proposition}
 Let $(C,G)$ be a big action such that  $G_2 \simeq (\Z/p\Z)^n$, with $n \geq 2$.
 We keep the notation defined above. 
 \begin{enumerate}
 \item For all $i$ in $\{1,\cdots, n-1\}$, $\ell_{i,i+1}$ is a linear form from $V$ to
 $\F_p$.
 \item For all $i$ in $\{1,\cdots, n-1\}$, put $\mathcal{L}_{i,i+1}(X):=\prod_{y \in Ker \ell_{i,i+1}} \,(X-y)$. Then, whenever $\ell_{i,i+1}$ is non identically zero, there exists $\lambda_i$ in $k-\{0\}$ such that, for all $y$ in $V$, $\ell_{i,i+1}(y)= \lambda_i \, \mathcal{L}_{i,i+1}(y)$. In this case,  $V=Z(\lambda_i^p\, \mathcal{L}_{i,i+1}^p-\lambda_i\, \mathcal{L}_{i,i+1})$.
 \end{enumerate}
 \end{proposition}

\noindent \textbf{Proof:} 
The matricial multiplication first ensures that for all $i$ in $\{1,\cdots, n-1\}$, $\ell_{i,i+1}$ is a linear form from $V$ to $\F_p$.
Besides, from the preliminary remarks of section 1, we infer that $P_V(X):=\prod_{y\in V} \, (X-y)$ is a separable additive polynomial of degree $p^v$, where $v$ denotes the dimension of the $\F_p$-vector space $V$. Then, for all $i$ in $\{1,\cdots,n-1\}$, $\mathcal{L}_{i,i+1}(X):=\prod_{y \in Ker \ell_{i,i+1}} \,(X-y)$ is an additive polynomial which divides $P_V(X)$. We now assume that $\ell_{i,i+1}$ is a nonzero linear form. In this case, $\mathcal{L}_{i,i+1}(X)$ has degree $p^{v-1}$ and there exists $\lambda_i$ in $k-\{0\}$ such that for all $y$ in $V$, $\ell_{i,i+1}(y)= \lambda_i \, \mathcal{L}_{i,i+1}(y)$. Since for all $y$ in $V$, $\ell_{i,i+1}(y)$ lies in $\F_p$, then $\lambda_i^p \, \mathcal{L}_{i,i+1}^p- \lambda_i\, \mathcal{L}_{i,i+1}=\lambda_i^p \, P_V$. The claim follows. $\square$

\begin{remarque}
By duality with respect to the Artin-Schreier pairing, the adapted basis of $A$ fixed in Definition 2.3.4 gives a basis of $G_2$, say $\{g_1,\cdots,g_n\}$, in which, the matrix of the automorphism $\phi(y)$ is the transpose matrix of $L(y)$ for all $y$ in $V$, namely a lower triangular matrix of 
$Gl_n(\F_p)$ with identity on the diagonal.
\end{remarque}

\begin{proposition}
Let $(C,G)$ be a big action such that  $G_2 \simeq (\Z/p\Z)^n$, with $n \geq 1$.
 We keep the notation defined above. 
For all integer $d$ such that $1 \leq d \leq n$, we denote by $A_d$ the $\F_p$-subvector space of $A$ generated by $\{\overline{f_i(X)}, \, 1 \leq i \leq d\}$. Let $H_d$ be the orthogonal of $A_d$ with respect to the Artin-Schreier pairing, namely the $\F_p$-subvector space of $G_2$ spanned by
$\{g_i,\, d+1 \leq i \leq n\}$ if $d<n$ and $H_n=\{0\}$. Then, the pair $(C/H_d, G/H_d)$ is a big action such that $(\frac{G}{H_d})_2=\frac{G_2}{H_d}$. It follows that $|\frac{G}{G_2}|=|\frac{G/H_d}{(G/H_d)_2}|$ and that the exact sequence
$$ 0 \longrightarrow G_2  \longrightarrow G  \stackrel{\pi}{\longrightarrow} V \longrightarrow 0$$
induces the following one:
$$ 0 \longrightarrow (G/H_d)_2 \simeq (\Z/p\, \Z)^d \longrightarrow G/H_d  \stackrel{\pi}{\longrightarrow} V \longrightarrow 0$$
\end{proposition}

\noindent \textbf{Proof:}
Since $\rho(V) \subset
T_{1,n}^u(\F_p)$, $A_d$ is stable under the action of $\rho$, that is to say under the translation: $X \rightarrow X+y$, with $y \in V$. By duality, $H_d$ is stable under the action of $\phi$, i.e. by conjugation by the elements of $G_1$. It follows that $H_d$ is a subgroup of $G_2$, normal in $G_1$. In this case, \cite{MR} (see Lemma 2.4 and Theorem 2.6) implies that the pair $(C/H_d, G/H_d)$ is a big action with $(\frac{G}{H_d})_2=\frac{G_2}{H_d} \subset \frac{G}{H_d}$. The claim follows. $\square$

\begin{corollaire}
Let $(C,G)$ be a big action such that  $G_2 \simeq (\Z/p\Z)^n$, with $n \geq 1$. Let the functions $f_i(X) \in k[X]$ be as in section 2.3.1.
Then, $f_1(X)=X \,S_1(X)+c\, X$, where $S_1 \in k\{F\}$ is an additive polynomial.
Furthermore, after an homothety and a translation on $X$, one can assume that $S_1$ is monic and $c=0$.
\end{corollaire}

\noindent \textbf{Proof:}
The function field of the curve $C/H_d$, as defined in Proposition 2.11, is parametrized by the $d$ Artin-
Schreier equations: $W_i^p-W_i=f_i(X)$, with $1 \leq i \leq d$.
 In particular, for $d=1$ $(C/H_1, G/H_1)$ is a big action whose second lower ramification group has
 order $p$. Then, \cite{LM} asserts that $f_1(X)=X \,S_1(X)+c\, X$ in $k[X]$, where $S_1 \in k\{F\}$
is an additive polynomial. $\square$

\subsection{Characterization of the trivial representation.}
\indent To conclude this section, we give a characterization of the case where the representation $\rho$ or $\phi$ is trivial.

\begin{proposition} 
Let $(C,G)$ be a big action such that $G_2\simeq (\Z/p\, \Z)^n$, with $n \geq 1$.\\  When
keeping the notation defined above, the following assertions are equivalent.
\begin{enumerate}
\item The representation $\phi$ is trivial, namely $\phi(V)=\{\id \}$.
\item The second ramification group $G_2$ is included in  the center of $G_1$.
\item The representation $\rho$ is trivial, i.e. $$\forall \, i \in \{1,\cdots,n\}, \quad \forall \, y \in V, \quad f_i(X+y)-f_i(X)=0 \mod \; \wp(k[X])$$
\item For all $i$ in $\{1,\cdots,n\}$, the functions $f_i$ read: $f_i(X)=X \, S_i(X)+ c_i\,X$ mod $\wp (k[X])$, where $S_i$ is an additive polynomial of degree $s_i\geq 1$ in $F$. Write $ S_i(F) = \sum _{j=0}^{s_i} \,a_{i,j} \,F^j$ with $a_{i,s_i} \neq 0$. Then, following \cite{El} (section 4), we can define an additive polynomial related to $f_i$, called the "palindromic polynomial" of $f$:
 $$Ad_{f_i}:= \frac{1}{a_{i,s_i}} \; F^{s_i}  (\sum _{j=0}^ {s_i} \,a_{i,j} \,F^j + F^{-j} \,a_{i,j}).$$
In this case, $$V \subset \bigcap_{i=1}^n Z (Ad_{f_i}) $$
\end{enumerate}
\end{proposition}

\noindent The proof of this proposition requires a preliminary lemma.

\begin{lemme}
When keeping the notation defined above, $\cap_{y \in V} Ker\, ( \phi(y)-\id) =Z(G_1) \cap G_2$.
\end{lemme}

\noindent \textbf{Proof of Lemma 2.14:} 
Consider $g$ in $G_2$. Then, $g$ lies in  $\cap_{y \in V} Ker\, ( \phi(y)-\id)$ if and only if
$\phi(y)(g)=g$ for all $y$ in $V$. For all $g_1$ in $G$, put $y_1:=\pi(g_1)$. By definition, the equality $\phi(y_1)(g)=g$ means that $g_1^{-1}\,g\, g_1=g$. This proves the expected formula. $\square$.

\medskip

\noindent \textbf{Proof of Proposition 2.13:}
The equivalence between the first and the second assertion derives from Lemma 2.14. As the equivalence between the first and the third point comes from the duality of $\phi$ and $\rho$ (cf. Proposition 2.2), the only point that has to be explained is the equivalence between the last assertion and the three preceding ones. \\
\indent For all $i$ in $\{1,\cdots,n\}$, the function field of the curve $W_i^p-W_i=f_i(X)$
is a $p$-cyclic \'etale cover of the affine line, whose Galois group is denoted by $H_i$. Then,
$H_i$ is an index $p$-subgroup of $G_2$. Besides,
if the second point is satisfied, $G_2$ is included in $Z(G_1)$, which implies that $H_i$ is normal in $G_1$. From \cite{MR} (see Lemma 2.4 and Theorem 2.6), we infer that $(C/H_i, G/H_i)$ is a big action whose second ramification group $(G/H_i)_2=G_2/H_i$ is $p$-cyclic. By \cite{LM} (see Prop. 8.3), $f_i(X)=c_iX +X S_i(X)$ mod $\wp(k[X])$, with $S_i \in k\{F\}$. In addition, $V$ is included in $Z(Ad_{f_i}).$ 
Conversely,  if $f_i(X)=X \;S_i(X)+c_iX $, then it follows from Proposition 5.5 in \cite{LM} that $Z(Ad_{f_i})= \{y \in k, \, f_i(X+y)-f_i(X)=0 \, \mod \, \wp(k[X]) \, \}$. Thus, the third point is verified. $\square$

\section{The link with the $k$-algebra generated by additive polynomials.}

\indent The purpose of this section is to highlight the role played by the $k$-algebra generated by the additive polynomials in the parametrization of big actions with a $p$-elementary abelian second ramification group.

\subsection{The $k$-algebra generated by additive polynomials.}

\begin{defi}
We define $\Sigma_1$ as the $k$-subvector space of $k[X]$ generated by $1$ and by the additive polynomials of $k[X]$.
More generally, for any $n \geq 1$, we define $\Sigma_n$ as the $k$-subvector space of $k[X]$ generated by $1$ and the products of at most $n$ additive polynomials of $k[X]$. 
For $n=0$, we put $\Sigma_0=k$ and for $n<0$, we put $\Sigma_n=\{0\}$.
\end{defi}

\begin{remarque}
\begin{enumerate}
\item For $n \geq 1$, this definition means that $f$ is a polynomial of $\Sigma_n$ if and only if there is a way to write $f$ as a linear combination over $k$ of products of at most $n$ additive polynomials. 
\item The sequence $(\Sigma_n)_{n \in \Z}$ enjoys the following properties:
 \begin{enumerate}
\item $1 \in \Sigma_0$ 
\item For all integer $n$ in $\Z$, $\Sigma_n \subset \Sigma_{n+1}$
\item For all integer $m$ and $n$ in $\Z$, $\Sigma_m \, \Sigma_n \subset \Sigma_{m+n}$.
\item $\bigcup_{n\in \Z} \Sigma_i=k[X]$
\end{enumerate}
 In particular, the sequence $(\Sigma_n)_{n \in \Z}$ is an increasing ring filtration of $k [X]$.
\end{enumerate}
\end{remarque}

 For a given $f$ in $k[X]$, we search for the minimal integer $n$ such that $f$ belongs to $\Sigma_n$. It requires
 the introduction of the order function related to the ring filtration.
 
\begin{defi}
Let $a$ be an integer whose $p$-adic expansion reads: 
$a= a_0+a_1\,p+a_2\,p+\cdots+a_t\,p^t$,
with $t \in \N$ and $0\leq a_i \leq p-1$, for all $i \in \{0,1,2,\cdots, t\}$. 
We define the integer $S_p(a) \in \N$ as the sum of the digits of $a$, namely:
$$S_p(a):=a_0+a_1+a_2+\cdots +a_t.$$
\end{defi}

\begin{remarque}
For all integer $m$ in $\N$, $S_p(m) = (p-1)\, v_p(m!)$, where $v_p$ denotes the $p$-adic valuation. 
We gather that, if $m_1$ and $m_2$ are two non-negative integers , $S_p(m_1+m_2) \leq S_p(m_1)+ S_p(m_2)$. 
\end{remarque}

\begin{lemme}
We keep the same notation as above. Let $a \in \N$ and $n \in \N$. 
Then, the monomial $X^a$ lies in $\Sigma_n$ if and only if $S_p(a) \leq n$.
It follows that $\inf \{n \in \N, \, X^a \in \Sigma_n\}= S_p(a)$.
\end{lemme}

\noindent \textbf{Proof:}
Assume that $X^a \in \Sigma_n$. It means that $X^a$ is a linear combination over $k$ of monomials of the form $X^{p^{\gamma_1} +p^{\gamma_2}+\cdots + p^{\gamma_t}}$, with $t \leq n$ and $\gamma_1 \leq \gamma_2 \leq \cdots \leq \gamma_t$. It follows that $a$ also reads $a=p^{\alpha_1} +p^{\alpha_2}+\cdots + p^{\alpha_t}$ with $t \leq n$ and $\alpha_1 \leq \alpha_2 \leq \cdots \leq \alpha_t$. Therefore, Remark 3.4 implies $S_p(a) = S_p(p^{\alpha_1} +p^{\alpha_2}+\cdots + p^{\alpha_t})\leq  S_p(p^{\alpha_1}) +S_p(p^{\alpha_2})+\cdots + S_p(p^{\alpha_t}) = t \leq n$.\\
\indent Conversely, we suppose that $S_p(a) \leq n$ and prove the result by induction on $n$. 
If $n=0$, then $S_p(a)=0$ and $X^a=X^0=1 \in \Sigma_0$.
We now assume that the property is true for $n$ and suppose that $S_p(a) \leq n+1$. If $S_p(a)=n$, then, by induction hypothesis, $X^a \in \Sigma_n \subset \Sigma_{n+1}$. Otherwise, $S_p(a)=n+1$ and there exists an integer $a_i$ in the $p$-adic expansion of $a$ such that $a_i \geq 1$. Put $b:=a-p^i$. As $S_p(b)=n$, the hypothesis implies $X^b \in \Sigma_n$, hence $X^a=X^b \, X^{p^i} \in \Sigma_{n+1}$. $\square$
\medskip

\begin{defi} Let $f$ be a nonnull polynomial of $k[X]$ such that $f(X)= \sum_{a \in \N} c_a(f)\, X^a$.
We define $$d_p(f):= \max_{c_a(f)  \neq 0} \, \{S_p(a) \}$$
By convenience, put $d_p(0):=-\infty$.
\end{defi}

\begin{lemme}
Let $f$ and $g$ be polynomials of $k[X]$. Let $n \in \Z$. We keep the same notation as above. 
\begin{enumerate}
\item  $f \in \Sigma_n$ if and only if $d_p(f) \leq n$.
In other words,  $f(X)= \sum_{a \in \N} c_a(f)\, X^a \in \Sigma_n$ if and only if,
whenever $S_p(a)>n$, $c_a(f)=0$.
\item If $f$ is non identically zero, $d_p(f)=\inf \{n \in \Z, \, f \in \Sigma_n\}$.
\item $d_p(f)=-\infty$ if and only if $f \in \cap_{n \in \Z} \Sigma_n=\{0\}$.
\item $d_p(f\,g) \leq d_p(f)+d_p(g)$. 
\item $d_p(f+g) \leq \sup \{ d_p(f), \,d_p(g)\, \}$. 
\item $d_p(F(f))= d_p(f)$, where $F$ means the Frobenius operator.
\item Let $S(X) \in k[X]$ be an additive polynomial. Then, $d_p(f(S(X)))=d_p(f(X))$.
\end{enumerate}
In particular, $d_p$ is the order function of the ring filtration defined by the $(\Sigma_n)_{n \in \Z}$.
\end{lemme}

\noindent \textbf{Proof:} Most of the properties can be deduced from Remark 3.4 and Lemma 3.5. The last one is left as an exercise to the reader. $\square$

\begin{defi} Let $f$ be a polynomial of $k[X]$. Let $y \in k$. We define the operator $\Delta_y$ as follows:
$\Delta_y(f):=f(X+y)-f(X)$.
\end{defi}

One checks that this operator enjoys the following property: 
\begin{lemme}
For all $y$ in $k$ and for all $n \in \Z$, $\Delta_y(\Sigma_{n+1}) \subset \Sigma_n$.
\end{lemme}

\begin{remarque} Although $d_p(\Delta_y(X^a))=d_p(X^a)-1$, for all $y$ in $k-\{0\}$ and all $a$ in $\N^*$,
one can find some polynomial $f$ in $k[X]$ and some $y$ in $k-\{0\}$ such that $d_p(\Delta_y(f))\neq d_p(f)-1$. It means that for $n \geq 2$ and for $y$ in $k-\{0\}$, $\Delta_y(\Sigma_{n+1}-\Sigma_n)$
is not always included in $\Sigma_n-\Sigma_{n-1}$.
\end{remarque}

\subsection{Notation and preliminary lemmas.}
\indent We begin by recalling some notation and proving some lemmas useful for the proof of next theorem. Let $(C,G)$ be a big action such that $G_2 \simeq (\Z/p\Z)^n$, with $n\geq 1$.
We call condition $(N)$ the inequality satisfied by big actions, namely: 
$\frac{|G|}{g} > \frac{2\,p}{p-1}$.
We fix an adapted basis of $A$: $\{\overline{f_1(X)}, \cdots, \overline{f_n(X)}\}$,
 as constructed in Definition 2.3 and assume that the functions $f_i(X)$ are reduced mod
 $\wp(k[X])$ (see definition in section 1). We denote by $m_i$ the degree of $f_i(X)$. As recalled in Corollary 2.12, $f_1(X)=X\, S_1(X)+c_1\,X$, where $S_1\in k\{F\}$ is an additive polynomial with degree
 $s_1 \geq 1$ in $F$. In this case, the palindromic polynomial $Ad_{f_1}$ related to $f_1$ is defined as in Proposition 2.13. Besides, the function field $L:=k(C)$ is parametrized by the $n$ Artin-Schreier equations: $W_i^p-W_i=f_i(X)$, with $1 \leq i \leq n$. We denote by $\rho$ the representation from $V$ to $Aut(A)$ defined in section 2.2. In the adapted basis
fixed above, the automorphism $\rho(y)$ is associated with the unipotent matrix:
$$ L(y):=
\begin{pmatrix}
 1 &\ell_{1,2}(y)&  \ell_{1,3}(y) &\cdots& \ell_{1,n}(y) \\
 0 & 1 &\ell_{2,3}(y) & \cdots&\ell_{2,n}(y) \\
 0 & 0 &\cdots&\cdots&\ell_{i,n}(y) \\
 0 & 0 &\cdots& 1& \ell_{n-1,n}(y)\\
 0 & 0 & 0 & 0 & 1
\end{pmatrix} \in Gl_n(\F_p)
$$

\begin{lemme} We keep the notation defined above.
The dimension of the $\F_p$-vector space $V$ satisfies $v \leq 2\,s_1$ and $p^v \geq m_n+1$.
In particualr, $2 \leq s_1+1 \leq v \leq 2\,s_1$.
\end{lemme}

\noindent \textbf{Proof:} The inclusion of $V$ in $Z(Ad_{f_1})$ first requires $v \leq 2\,s_1$. On the one hand, $|G|=p^{n+v}$. On the other hand, Corollary 2.7 implies: $g=\frac{p-1}{2} \, \sum_{i=1}^n \, p^{i-1}\, (m_i-1) \geq \frac{p-1}{2}\, p^{n-1} \, (m_n-1)$. Thus, $\frac{|G|}{g} \leq \frac{2\,p}{p-1}\, \frac{p^v}{m_n-1}.$  The inequality $p^v \leq m_n-1$  would contradict condition $(N)$. Therefore, since $m_n$ is prime to $p$, we obtain $p^v \geq m_n+1$.
It follows that $p^v > m_n \geq m_1=1+p^{s_1}$ and $v \geq s_1+1$. $\square$
\medskip

\begin{lemme}
Let $f(X):=\sum_{a \in \N} c_a(f)\, X^a$ be a polynomial in $ \wp(k[X])$. Fix $a_0 \in \N-p\,\N$ and define $I_{a_0}:=\{a_0\,p^n, \,n \in \N\}$.
Then, the polynomial $f_{a_0}(X):= \sum_{a \in I_{a_0}} c_a(f)\, X^a$ also lies in $\wp(k[X])$. In particular, if $f_{a_0}(X)$ is non identically zero, then $p$ divides its degree.
\end{lemme}

\noindent \textbf{Proof:}
The Frobenius operator $F$ acts on the basis $(X^a)_{a\in \N}$ of $k[X]$ and this action induces a partition
of the monomials of $k[X]$, namely $(X^a)_{a \in I_{a_0}}$, for $a_0$ running over $\{0\} \cup  \{\N-p\,\N\}$. This justifies the first claim.
Now, assume that $f_{a_0}(X)$ is non identically zero. If $f= \wp(g)$ with $g \in k[X]$, then $f_{a_0}= \wp(g_{a_0})$, with
$g_{a_0}$ defined as for $f$. It follows that $deg \, (f_{a_0})=p \, deg \, (g_{a_0})$. $\square$

\subsection{The link with the parametrization of big actions.}

\begin{theoreme}
We keep the notation defined in sections 3.1 and 3.2.\\
For all $i$ in $\{1,\cdots,n\}$, $f_i(X)$ belongs to $\Sigma_{i+1}$.
\end{theoreme}

\noindent \textbf{Proof:}
For a fixed $n$, we proceed by induction on $i$. 
As recalled above, $f_1(X)=X S_1(X)+c_1\, X$, where $S_1$ is an additive polynomial. Accordingly, $f_1 \in \Sigma_2$. We now consider some integer $i$ such that $ 2 \leq i \leq n$ and assume that for all $j$ in $\{1,\cdots,i-1\}$, $f_j(X)$ lies in $\Sigma_{j+1}$. From the form of the matrix $L(y)$, we gather:
$$
\forall \, y \in V, \; \Delta_y(f_i):=f_i(X+y)-f_i(X)=\sum_{j=1}^{i-1} \,\ell_{j,i}(y) \,f_j(X) \quad 
\mod \,(\wp(k[X]))$$
where for all $j$ in $\{1,\cdots,i-1\}$, $\ell_{j,i}$ is a map from $V$ to $\F_p$. \\
\indent Suppose that $f_i(X)$ does not belong to $\Sigma_{i+1}$ and call $X^a$ the monomial of $f_i(X)$ with highest degree which does not belong to $\Sigma_{i+1}$. Note that, by definition of $a$, $a \geq i+1$. Furthermore, as $f_i$ is assumed to be reduced mod $\wp(k[X])$, $a \neq 0$ mod $p$.\\
\indent We first prove that $p$ divides $a-1$
Indeed, assume that $p$ does not divide $a-1$ and apply Lemma 3.12 to $f(X):=\Delta_y(f_i)-\sum_{j=1}^{i-1} \,\ell_{j,i}(y) \,f_j(X)$ and $a_0:=a-1 \in \N-p\N$. 
To construct the polynomial $f_{a_0}$ as defined in Lemma 3.12, we first search for the monomials $X^{(a-1)p^r}$, with $r \geq 0$, in $\Delta_y(f_i)$. 
If $r>0$, such monomials come from monomials $X^b$ of $f_i(X)$ such that $b >(a-1)\,p^r \geq (a-1)\,p
 \geq a$, since $a \geq i+1 \geq 2 \geq \frac{p}{p-1}$.
By definition of $a$, such monomials $X^b$, whose degree is strictly higher than $a$, lies in $\Sigma_{i+1}$. Then, by Lemma 3.9, they generate in $\Delta_y (f_i)$ polynomials which belongs to $\Sigma_i$. But $X^{a-1} \not \in \Sigma_i$: otherwise, $X^a \in X \Sigma_i \subset \Sigma_{i+1}$, which contradicts the definition of $a$. We infer from Lemma 3.7.6 that no $X^{(a-1)p^r}$, with $r \geq 0$, lies in $\Sigma_{i}$.
It follows that no monomial $X^{(a-1)p^r}$, with $r>0$, can be found in $\Delta_y(f_i)$. We now search for the monomial $X^{a-1}$. By the same token, one can check that the only monomial in $f_i(X)$ which generates $X^{a-1}$ in $\Delta_y(f_i)$ is $X^a$. More precisely, it produces $a \, y \, c_a(f_i)\, X^{a-1}$ in $\Delta_y(f_i)$, where $c_a(f_i) \neq 0$ denotes the coefficient of $X^a$ in $f_i$. 
As the induction hypothesis asserts that $\sum_{j=1}^{i-1} \,\ell_{j,i}(y) \,f_j(X)$ lies in $\Sigma_i$, which is the case of none of the $X^{(a-1)p^r}$, we gather that $f_{a_0}(X)=a\,y \,c_a(f_i) \, X^{a-1}$. As $p$ does not divide $a_0=a-1$, it follows from Lemma 3.12 that $f_{a_0}(X)$ is identically zero. Since $a \neq 0$ mod $p$, this implies that $y=0$ for all $y$ in $V$, hence $V=\{0\}$. It means that $G_1=G_2$, which is impossible for a big action. Accordingly, $p$ divides $a-1$. 
Thus, we can write $a=1+\lambda \; p^t$ with $t \geq 1$, $\lambda$ prime to $p$ and $\lambda >i \geq 2$, as $X^a$ does not lie in $\Sigma_{i+1}$. \\
\indent Now, put $j_0:=a-p^t=1+(\lambda-1)\, p^t$ and apply Lemma 3.12 to $f(X):=\Delta_y(f_i)-\sum_{j=1}^{i-1} \,\ell_{j,i}(y) \,f_j(X)$ and $a_0:=j_0 \in \N-p\,\N$.
To construct the polynomial $f_{a_0}$, we first determine the monomials $X^{j_0 p^r}$, with $r \geq 0$, occuring in $\Delta_y(f_i)$. 
If $r>0$, such terms come from monomials $X^b$ of $f_i(X)$ such that $b >j_0p$. 
But $j_0p > a$. Indeed, $$j_0p \leq a \Leftrightarrow p\,(1+(\lambda-1)\,p^t) \leq 1+ \lambda\,p^t \Leftrightarrow \lambda \leq \frac{1-p+p^{t+1}}{p^t\,(p-1)}=\frac{-1}{p^t}+ \frac{p}{p-1} < \frac{p}{p-1} \leq 2$$ which contradicts $\lambda \geq 2$.
As explained above, the monomials $X^b$ of $f_i(X)$, with $b >a$, produce polynomials in $\Delta_y(f_i)$ which belongs to $\Sigma_{i}$, whereas $X^{j_0}$ does not belong to $\Sigma_i$. Otherwise, $X^a=X^{p^t} X^{j_0}$ would belong to $\Sigma_{i+1}$, hence a contradiction. We gather from Lemma 3.7.6 that no $X^{j_0 p^r}$, with $r \geq 0$, lies in $\Sigma_{i+1}$.
It follows that no monomials $X^{j_0 p^r}$, with $r>0$, can be found in $\Delta_y(f_i)$.
Likewise, for $r=0$, the only monomials of $f_i(X)$ which generates $X^{j_0}$ in $\Delta_y(f_i)$ are those of the form: $X^b$, with $j_0+1  \leq b \leq a$. For all $b \in \{j_0+1,\cdots, a\}$, the monomial $X^b$ of $f_i(X)$ generates some $\binom b{j_0}\, y^{b-j_0}\, X^{j_0}$ in $\Delta_y(f_i)$. 
It follows that the coefficient of $X^{j_0}$ in $\Delta_y(f_i)$ is $T(y)$ with $T(Y):= \sum _{b=j_0+1}^a  \, c_b(f_i) \,\binom b{j_0} \,Y^{b-j_0}$, where $c_b(f_i)$ denotes the coefficient of $X^b$ in $f_i(X)$.  As no $X^{j_0p^r}$, with $r\geq 0$, can be found in $\sum_{j=1}^{i-1} \,\ell_{j,i}(y) \,f_j(X)$ which lies in $\Sigma_i$ by induction, the polynomial $f_{a_0}$ eventually reads $f_{a_0}(X)=T(y)\,X^{a_0}$. By Lemma 3.12, $f_{a_0}$ is identically zero, which means that for all $y$ in $V$, $T(y)=0$. We gather that $V$ is included in the set of zeroes of $T$. As the coefficient of $Y^{a-j_0}$ in $T(Y)$ is $c_a(f_i) \binom a{j_0}=c_a(f_i)
\binom {1+\lambda p^t}{1+(\lambda-1)p^t} \equiv c_a(f_i) \lambda \neq 0$ mod $p$, the polynomial $T(Y)$ has degree $a-j_0=p^t$,  hence $v\leq t$.
This leads to a contradiction, insofar as Lemma 3.11 implies:
$p^v \geq m_n-1 \geq m_i-1 \geq a-1 = \lambda \, p^t \geq 2 \, p^t >p^t $, which involves: $v >t$. 
As a consequence, $f_i(X)$ does not have any monomial which does not belong to $\Sigma_{i+1}$, which completes both the induction and the proof of the theorem.
$\square$

\begin{remarque}
The proof is actually self-contained, since the first step of the induction, namely $f_1 \in \Sigma_2$, could be obtained without any hint at Corollary 2.12 which requires the use of \cite{LM}.
Indeed, in the case $i=1$, the sum $\sum_{j=1}^{i-1} \,\ell_{j,i}(y) \,f_j(X)$ is replaced by $0$ which
obviously lies in $\Sigma_1$. Using the same argument as in the second part of the proof, it enables us to conclude that $f_1$ belongs to $\Sigma_2$.
\end{remarque}

\section{A special case of trivial representation.}
\indent This section is devoted to a first special case in which the representation $\rho$ is trivial or, equivalently, each function $f_i$ lies in $\Sigma_2-\Sigma_1$. 
The difficulty in solving the general case of trivial representation lies in finding the GCD for the family of
palindromic polynomials associated to the function $f_i$ as defined in Proposition 2.13. This could be done by working in the Ore ring of Laurent polynomials $k\{F,F^{-1}\}$ (see \cite{El}, section 3, or \cite{Go}, 1.6). Nevertheless, in what follows, we merely explore the simplest case
where all the palindromic polynomials are equal.

\begin{lemme}
Let $(C,G)$ be a big action such that $G_2 \simeq (\Z/p\Z)^n$, with $n \geq 1$.
We keep the notation defined in sections 2.3.2 and 4.2.
Then, the following assertions are equivalent.
\begin{enumerate}
\item The upper ramification filtration of $G_2$ has only one jump.
\item The functions $f_i$'s have the same degree, i.e. for all $i$ in $\{1,\cdots,n\}$, $m_i=m_1=1+p^{s_1}$.
\end{enumerate}
In this case, the representation $\rho$ is trivial and each function $f_i$ reads $f_i(X)=X\,
 S_i(X)+c_i\,X \, \in \Sigma_2-\Sigma_1$, where $S_i$ is an additive polynomial with degree $s_1$ in $F$. Moreover, $V \subset \cap_{1\leq i
\leq n} \, Z(Ad_{f_i})$. 
\end{lemme}

\noindent \textbf{Proof:}
Assume that there is only one jump in the upper ramification filtration of $G_2$ as defined in section 2.3.2, namely $G_2=(G_2)^{\nu_0} \supsetneq (G_2)^{\nu_0+1}=\{0\}$. 
The duality between the filtrations of $A$ and $G_2$ (cf. Proposition 2.5) implies that this is equivalent to $\{0\}=A^{\mu_0}\subsetneq A^{\mu_0+1}=A$. By Remark 2.4, this situation occurs if and only if all the functions $f_i(X)$ have the  same degree, namely: $1+p^{s_1}$. In this case, it follows from Remark 2.8 that the representation $\rho$ is trivial. Then, the following assertions derive from Proposition 2.13. $\square$  
\smallskip

 In what follows, we restrict to the special case: $V=Z(Ad_{f_1})$, which means that $V$ has maximal cardinality for a given $s_1$, namely $|V|=p^{2s_1}$.

\begin{proposition} Let $n \geq 2$.
 We assume that $\rho(V)=\{\id\}$ and keep the notation defined above. We suppose that $v=2\,s_n$.
 Then, for all $i$ in $\{1,\cdots,n\}$, $s_i=s$ and $V=Z(Ad_{f_1})$.
 Furthermore, there exists an integer $d$ dividing $s$ and some $\gamma_2, \cdots, \gamma_n$ in
 $\F_{p^d}-\F_p$ such that: $$S_1= \sum_{j=0}^{s/d} \, a_{jd} \,F^{jd}
 \qquad and \qquad \forall \, i \in \{2,\cdots,n\} \quad S_i= \gamma_i\; S_1$$
Moreover, $\{ \gamma_1:=1, \gamma_2, \cdots, \gamma_n\}$ are linearly independent over $\F_p$.
It follows that $s \geq 2$.
\end{proposition}

\noindent \textbf{Proof:}
As $v \leq 2\,s_1 \leq 2\,s_n$, the hypothesis $v=2\,s_n$ implies that each $s_i$ is equal to $s_1$. From now on, $s_1=s_2=\cdots=s_n$ is denoted by $s$. By Proposition 2.13, $V \subset \cap_{1\leq i \leq n} Z(Ad_{f_i}) \subset Z(Ad_{f_1})$. As the two vector spaces $V$ and $Z(Ad_{f_1})$ have the same dimension over $\F_p$, namely $v=2\,s$, we conclude that $Z(Ad_{f_1})=V=Z(Ad_{f_i})$ for all $i$ in $\{1,\cdots,n\}$. Since $k$ is algebraically closed and since $Ad_{f_1}$ and $Ad_{f_i}$ are monic, it follows that $Ad_{f_1}=Ad_{f_i}$.\\
\indent Let $i$ in $\{2,\cdots, n\}$. Write: $S_1= \sum_{k=0}^s a_k\; F^k$ and $S_i= \sum_{k=0}^s b_k\; F^k$, with $a_s \neq 0$ and $b_s \neq 0$.
Then, $Ad_{f_1} = \frac{1}{a_s}\;  F^s \; (\sum_{k=0}^s (a_k \; F^k + F^{-k} \; a_k))= Ad_{f_i} = \frac{1}{b_s}\;  F^s \; (\sum_{k=0}^s (b_k \; F^k + F^{-k} \; b_k))$. As for all $\alpha \in k$, $F \; \alpha=\alpha^p \; F$, we obtain:
$ \gamma_i\; \sum_{k=0}^s (a_k^{p^s} \; F^{s+k} + a_k^{p^{s-k}} \; F^{s-k}) = \sum_{k=0}^s (b_k^{p^s} \; F^{s+k} + b_k^{p^{s-k}} \; F^{s-k})$, with $\gamma_i\, a_k^{p^s}= b_k^{p^s}$ and $\gamma_i \,a_k^{p^{s-k}}= b_k^{p^{s-k}}$, for all $0 \leq k \leq s$. It implies that 
$\gamma_i^{p^k}=\gamma_i$, for all $0 \leq k \leq s$ such that $a_k \neq 0$. In particular, as $a_s \neq 0$, then $\gamma_i \in \F_{p^s}$. If we denote by $d$ the degree of the minimal polynomial of $\gamma_i$ over $\F_p$, then  $\F_{p^d}:=\F_p(\gamma_i) \subset \F_{p^s}$, so $d$ divides $s$.
By the same token, for all $0 \leq k \leq s$ such that $a_k \neq 0$, $\gamma_i \in \F_{p^k}$. Therefore, $\F_{p^d}=\F_p(\gamma_i) \subset \F_{p^k}$, which proves that $d$ divides $k$, whenever $a_k \neq 0$. It follows that $S_1= \sum_{j=0}^{s/d}\, a_{jd} \,F^{jd}$. In addition, as $\gamma_i \in \F_{p^s}$, we gather from $b_k^{p^s}=\gamma_i \,a_k^{p^s}$ for all $0 \leq k \leq s$, that $S_i =\gamma_i \, S_1$. \\
\indent Note that $\{\gamma_1,\cdots, \gamma_n\}$ are linearly independent over $\F_p$. Otherwise, it would contradict the condition (c) imposed in Definition 2.3.4. It follows that none of the $\gamma_i$'s, for $i\geq 2$, are in $\F_p$. So $s\geq 2$. $\square$

\medskip
We now display a family of big actions satisfying the conditions
described in Proposition 4.2.

\begin{proposition}
\begin{enumerate}
\item Let $s \in \N^*$, $d \in \N^*$ dividing $s$ and $n \in \N^*$ such that $n \leq d$. \\
Take $\{\gamma_1:=1, \gamma_2, \cdots, \gamma_n\}$ in $\F_{p^d}$, linearly independent over $\F_p$.\\
Put $S_1:=\sum_{j=0}^{s/d} \, a_{jd} \, F^{jd} \in k\{F\}$, with $a_s \neq 0$.
For all $i$ in $\{1,\cdots, n\}$, we define $S_i:=\gamma_i \, S_1$ in $k\{F\}$ and $f_i(X):=X\,S_i(X)+c_i \,X \in k[X]$. Then, for all $i$ in $\{1,\cdots, n\}$, $Z(Ad_{f_i})=Z(Ad_{f_1})$. Put $V:=Z(Ad_{f_1})$. 
\item The function field of the curve $C$ parametrized by: 
$W_i^p-W_i=f_i(X)$ with $1 \leq i \leq n$,
 is an \'etale extension of $k[X]$ with Galois group $H \simeq (\Z/p\Z)^n$. Then, the group of translations of the affine line: $\{\tau_y:X \rightarrow X+y, \, y \in V\}$ extends to an automorphism $p$-group of $C$, say $G$, such that:
$$0 \longrightarrow  H\simeq (\Z/p \, \Z)^n \longrightarrow G \longrightarrow V \longrightarrow 0$$
\item Thus, we obtain a big action $(C,G)$ whose second ramification group $G_2$ is isomorphic to $(\Z/ p\, \Z)^n$ and such that the representation
$\rho$ is trivial. 
 Moreover, $Z(G)=D(G)=Fratt(G)\simeq (\Z/p\Z)^n$, where $Fratt(G)$ denotes the Frattini
subgroup of $G$ and $D(G)$ its derived subgroup.
 Then, $G$ is a special group (see \cite{Su2}, Def. 4.14).
 Besides, if $p >2$, $G$ has exponent $p$, whereas, if $p=2$, $G$ has exponent $p^2$. 
\end{enumerate}
 \end{proposition}

\noindent \textbf{Proof:} 
\begin{enumerate}
\item Fix $i$ in $\{1,\cdots, n\}$. Then, $\gamma_i \,Ad_{f_1}= \frac{\gamma_i}{a_s} \, \sum_{j=0}^{s/d} \, (a_{jd} \, F^{jd} +F^{-jd} \, a_{jd})$.
Since $\gamma_i$ lies in $\F_{p^d}$,
$\gamma_i \,Ad_{f_1}= \frac{1}{a_s} \, \sum_{j=0}^{s/d} \, (a_{jd} \, \gamma_i \, F^{jd} +F^{-jd} \, a_{jd} \gamma_i)
= \frac{a_s\,\gamma_i}{a_s} \, Ad_{f_i}= \gamma_i \, Ad_{f_i}$. 
 So, $Z(Ad_{f_i})=Z(Ad_{f_1})$.  
\item As $Z(Ad_{f_i})=\{ y \in k, \, \Delta_y(f_i)=0 \, \mod
\wp(k[X]) \}$ (see \cite{LM}, Prop. 5.5), it follows that, for all $y$ in $V$, $\Delta_y(f_i)=0$ mod $\wp(k[X])$. So, Galois theory ensures the existence of the group $G$. 
 \item We deduce from the first point that $|G|=|G_2||V|=p^{n+2\,s}$. We compute the genus of $C$ by means of the formula given in Corollary 2.7. This yields: $g=\frac{(p^n-1)\, p^{s}}{2}$. Therefore, $\frac{|G|}{g}= \frac{2\,p^{n+s}}{p^n-1}$. So, the pair $(C,G)$ is a big action.\\
\indent We now show that $Z(G)=D(G)$. By Proposition 2.13, $Z(G)$ contains $G_2$ which coincides with $D(G)$ (see Theorem 2.6 in \cite{MR}).
 Conversely, let $\mathcal{H}$ be an index $p$-subgroup of $D(G)$. As $\mathcal{H}\subset D(G) \subset Z(G)$, $\mathcal{H}$ is normal in $G$ and  Lemma 2.4.2 in \cite{MR} implies that the pair $(C/\mathcal{H},G/\mathcal{H})$ is a big action. The curve $C/\mathcal{H}$ is
parametrized by an Artin-Schreier equation: $W^p-W=f(X):=\sum_{i=1}^n \, \lambda_i\,f_i(X)$, with $(\lambda_1,\cdots, \lambda_n) \in \F_p^n$
  not all zeros. Condition (c) of Definition 2.3.4 imposes $deg(f)= \max_{i=1,\cdots,n} \{deg \;  \lambda_i \,
  \,f_i(X)  \}=1+p^{s_1}$. Besides, by Proposition 2.11, we get the following exact sequence:
 $$ 0 \longrightarrow D(G /\mathcal{H}) \simeq \Z / p\Z \longrightarrow G/\mathcal{H}
 \longrightarrow V\simeq (\Z / p\Z)^{2s_1} \longrightarrow 0$$
 In this case, Proposition 8.1 in \cite{LM} shows that $G/\mathcal{H}$ is an extraspecial group, which involves 
 $D(G) /\mathcal{H}=D(G/\mathcal{H})=Z(G/\mathcal{H})$. We denote by $\pi: G \rightarrow
 G/\mathcal{H}$ the canonical mapping. Then, $\pi(Z(G)) \subset
 Z(G/\mathcal{H})=D(G) /\mathcal{H}$. As $\mathcal{H}$ is included in $Z(G)$, it follows that $Z(G)
 \subset D(G)$. Since $Z(G)=D(G)=G_2\simeq (\Z/p\Z)^n$ and since $D(G)=Fratt(G)$ for any $p$-group $G$, we gather that $G$ is a special group. Moreover, if $p>2$, Proposition 8.1 in \cite{LM} shows that $G/\mathcal{H}$ is an extraspecial group with exponent $p$. Then, $\pi(G)^p=\{e\}$ and $G^p$ is included in $\mathcal{H}$, for any hyperplanes $\mathcal{H}$ of $D(G)$. It follows that $G^p=\{e\}$. If $p=2$, the same proposition shows that $G/\mathcal{H}$ has exponent $p^2$. It implies that $G$ also has exponent $p^2$.
 $\square$
\end{enumerate}

\begin{remarque}
In the situation described in Proposition 4.3, $\frac{|G|}{g^2}= \frac{4\,p^{n}}{(p^n-1)^2}$. Note that the latter quotient does not depend on $s$ any more.
\end{remarque}

\begin{remarque}
For any big action $(C,G)$, the group $G$ is included in the wild inertia subgroup of $Aut_k(C)$ at $\infty$, denoted by $G_{\infty,1}$. Furthemore,
Corollary 2.10 in \cite{MR} shows that the pair $(C,G_{\infty,1} )$ is a big action with $D(G_{\infty,1})=D(G)$. Now, assume that $(C,G)$ is a big action as described in Proposition 4.3. Then,  $p^{2s}=\frac{|G|}{|D(G)|}\leq \frac{|G_{\infty,1}|}{|D(G_{\infty,1})|}\leq p^{2\,s}$.
It follows that, in this special case, $G$ is equal to $G_{\infty,1}$. 
\end{remarque}

\section{The special case: $f_i$ in $\Sigma_{i+1}-\Sigma_i$.}
\indent In this section, we define a filtration on the derived group $D(G)$ of any group $G$. Then, we investigate the case where $G$ is extension of two elementary abelian $p$-groups and where the number of jumps in this filtration is maximal. Knowing that, for a big action, $G_2=D(G)$ (see Theorem 2.6 in \cite{MR}), we apply these results to the case of big actions with a $p$-elementary abelian $G_2$. This allows us to give a group-theoretic condition to characterize the big actions such that each function $f_i$ lies in $\Sigma_{i+1}-\Sigma_i$. In this situation, we prove that the filtration on $D(G)$ actually coincides with the upper ramification filtration of $G_2$ as exposed in section 2.3.2. and that, as opposed to the previous case, the number of jumps in the filtration is maximal whereas the cardinality of $V$ is minimal in regard to Lemma 3.11, namely: $v=s_1+1$.

\subsection{A filtration on $D(G)$.}

\begin{defi}
For any group $G$, we define a sequence of subgroups $(\Lambda_i(G))_{i \geq 0}$ as follows.
Put $\Lambda_0(G):= \{e\}$, where $e$ means the identity element of $G$. For all $i \geq 1$,
let $\pi_{i-1}: G \rightarrow \frac{G}{\Lambda_{i-1}(G)}$ be the canonical mapping. Then, $\Lambda_i(G)$ is the subgroup of $G$ defined by $\pi_{i-1}^{-1} (Z( \frac{G}{\Lambda_{i-1}(G)}) \cap D( \frac{G}{\Lambda_{i-1}(G)}))$.
Therefore, $$ \frac{\Lambda_i(G)}{\Lambda_{i-1}(G)}=Z( \frac{G}{\Lambda_{i-1}(G)}) \cap D( \frac{G}{\Lambda_{i-1}(G)}).$$
In this way, we get an ascending sequence of subgroups of $D(G)$:
$$\{e\}= \Lambda_0(G) \subset \Lambda_1(G) \subset \Lambda_2(G) \subset \cdots \subset D(G)$$
which are characteristic subgroups of $G$.
\end{defi}

We study the filtation defined above in the special case where 
$G$ is a $p$-group with the exact sequence:
\begin{equation} \label{2}
 0 \longrightarrow D(G) \simeq (\Z/p \,\Z)^n \longrightarrow G \stackrel{\pi}{\longrightarrow} V \simeq (\Z/\,p\, \Z)^v \longrightarrow 0 
\end{equation}
In other words, $G$ is a $p$-group whose Frattini subgroup is equal to $D(G)\simeq (\Z/\,p\, \Z)^n$, with $n \geq 1$. For convenience, we fix a set theoritical section, i.e. a map $s: G/D(G) \rightarrow G$ such that
$\pi \circ s=id_{G/D(G)}$. We also define a representation 
$\phi: G/D(G) \rightarrow Aut(D(G))$ as follows. For all $y$ in $G/D(G)$
 and all $g$ in $D(G)$, $\phi(y)(g)=s(y)^{-1} \, g \, s(y)$. As $G/D(G)$ is a
 $p$-group, one can find a basis $\{g_1, \cdots, g_n\}$ of the $\F_p$-vector space $D(G)$ in which, for all $y$ in $G/D(G)$, the matrix of the automorphism $\phi(y)$ belongs to the subgroup of $Gl_n(\F_p)$ made of the lower triangular matrices with identity
 on the diagonal, namely:
  $$ \Phi(y):=
\begin{pmatrix}
 1 & 0 &  0 &\cdots& 0 \\
 \ell_{2,1}(y) & 1& 0 & \cdots  & 0 \\
 \ell_{3,1}(y)  & \ell_{3,2}(y) & \cdots&0 &0 \\
 \ell_{i,1}(y) & \ell_{i,2}(y) & \cdots & 1& 0\\
 \ell_{n,1}(y) & \ell_{n,2}(y) & \cdots & \ell_{n,n-1}(y) & 1
\end{pmatrix} \in Gl_n(\F_p)$$
Note that for $n \geq 2$ and for all $i$ in $\{1, \cdots, n-1\}$, $\ell_{i+1,i}$ is a linear form from $G/D(G)$ to $\F_p$.

\begin{proposition}
Let $G$ be a group satisfying \eqref{2}. We keep the notation defined above. Then, the following assertions are equivalent.
\begin{enumerate}
\item The filtration defined by the $(\Lambda_i)_{i \geq 0}$ satisfies:
 $$\{e\}= \Lambda_0(G) \subsetneq \Lambda_1(G) \subsetneq \Lambda_2(G) \subsetneq \cdots 
\subsetneq \Lambda_n= D(G)$$ which means, for all $i$ in $\{1,\cdots,n\}$,
$$ \frac{\Lambda_i(G)}{\Lambda_{i-1}(G)}=Z( \frac{G}{\Lambda_{i-1}(G)}) \cap D( \frac{G}{\Lambda_{i-1}(G)}) \simeq \Z/p\Z.$$
\item For all $i$ in $\{1,\cdots,n\}$, $\Lambda_i(G)$ is the $\F_p$-subvector space of 
$D(G)$ spanned by $\{g_{n-i+1}, \cdots, g_n\}$.
\item For $n \geq 2$ and for all $i$ in $\{1, \cdots, n-1\}$, $\ell_{i+1,i}$ is a nonzero linear form.
\end{enumerate}
\end{proposition}

\noindent \textbf{Proof:}
We prove that the first point implies the second one by induction on $i$. 
Assume $i=1$. By the same argument as in Lemma 2.14, one proves that $\Lambda_1(G)=Z(G) \cap D(G)$ is equal to $ \cap_{y \in G/D(G)} \, Ker (\phi(y)-\id)$.
Then, the form of $\Phi(y)$ shows that $\cap_{y \in G/D(G)} \, Ker (\phi(y)-\id)$ contains the $\F_p$-vector space spanned by $g_n$.
As $\Lambda_1(G)$ is assumed to be isomorphic to $\Z/p\Z$, it follows that $\Lambda_1(G)= \F_p\, g_n$.
Now, take $i \geq 2$ and assume that $\Lambda_{i-1}(G)$ is the $\F_p$-subvector space of 
$D(G)$ spanned by $\{g_{n-i+2}, \cdots, g_n\}$.
Then, $\frac{G}{\Lambda_{i-1}(G)}$ is a $p$-group with the following exact sequence:
$$ 0 \longrightarrow \frac{D(G)}{\Lambda_{i-1}(G)}=D(\frac{G}{\Lambda_{i-1}(G)}) \simeq (\Z/p \,\Z)^{n-i+1} \longrightarrow \frac{G}{\Lambda_{i-1}(G)} \stackrel{\pi}{\longrightarrow} G/D(G) \simeq (\Z/\,p\, \Z)^v \longrightarrow 0$$  
This exact sequence induces a representation $\phi_{i-1}: G/D(G) \rightarrow Aut(\frac{D(G)}{\Lambda_{i-1}(G)})$. Consider the canonical mapping: 
$\pi_{i-1}:D(G) \rightarrow \frac{D(G)}{\Lambda_{i-1}(G)}$. In the basis $\{\pi_{i-1}(g_1) \cdots, \pi_{i-1}(g_{n-i+1})\}$, the matrix $\Phi_{i-1}(y)$ of
 the automorphism $\phi_{i-1}(y)$ reads:
   $$ \Phi_{i-1}(y):=
\begin{pmatrix}
 1 & 0 &  0 &\cdots& 0 \\
 \ell_{2,1}(y) & 1& 0 & \cdots  & 0 \\
 \ell_{3,1}(y)  & \ell_{3,2}(y) & \cdots&0 &0 \\
 \ell_{i,1}(y) & \ell_{i,2}(y) & \cdots & 1& 0\\
 \ell_{n-i+1,1}(y) & \ell_{n-i+1,2}(y) & \cdots & \ell_{n-i+1,n-i}(y) & 1
\end{pmatrix} \in Gl_{n-i+1}(\F_p)$$ 
where the maps $\ell_{i,j}$ are the same as in $\Phi(y)$.
As in the case $i=1$, $\frac{\Lambda_i(G)}{\Lambda_{i-1}(G)}=Z( \frac{G}{\Lambda_{i-1}(G)}) \cap D( \frac{G}{\Lambda_{i-1}(G)})$ is equal to $\cap_{y \in G/D(G)}
Ker \, (\phi_{i-1}(y)-\id)$. The latter is the 
$\F_p$-vector space of $D(\frac{G}{\Lambda_{i-1}(G)})=\frac{D(G)}{\Lambda_{i-1}(G)}$ generated by $\pi_{i-1}(g_{n-i+1})$. It follows that 
$\Lambda_i(G)$ is the $\F_p$-subvector space of 
$D(G)$ spanned by $\{g_{n-i+1}, \cdots, g_n\}$. As the second assertion trivially implies the first one, the equivalence between 1 and 2 is established.\\
\indent We now prove that the second assertion implies the third one. 
Take $i \geq 1$. As seen above, $\frac{\Lambda_i(G)}{\Lambda_{i-1}(G)}=\cap_{y \in G/D(G)}
Ker \, (\phi_{i-1}(y)-\id)$ is the
$\F_p$-vector space spanned by $\pi_{i-1}(g_{n-i+1})$.
From the form of the matrix $\Phi_{i-1}(y)$, we gather that $\ell_{n-i+1,n-i}$ is non identically zero.
The proof of the converse works by induction on $i$. If $i=1$, the form of the matrix
 $\Phi(y)$, with each $\ell_{i+1,i}$ non identically zero, implies that
 $\Lambda_1(G)= \cap_{y \in G/D(G)} Ker \, (\phi(y)-\id)$ is the $\F_p$-subvector space
 of $D(G)$ spanned by $g_n$. Now, take $i \geq 2$ and assume that $\Lambda_{i-1}(G)$ is the
 $\F_p$-subvector space of $D(G)$ spanned by $\{g_{n-i+2}, \cdots, g_n\}$. By hypothesis, each linear form $\ell_{i+1,i}$ occuring in $\Phi_{i-1}(y)$
 is non
 identically zero. It implies that $\frac{\Lambda_i(G)}{\Lambda_{i-1}(G)} = \cap_{y \in G/D(G)} Ker \,
 (\phi_{i-1}(y)-\id)$ is the $\F_p$-vector space spanned by $\pi_{i-1}(g_{n-i+1})$. We conclude as above. $\square$

 \begin{remarque} Note that the third condition of the Proposition 5.2
  does not actually depend on the triangularization basis $\{g_1,\cdots,g_n\}$ chosen for $D(G)$.
 \end{remarque}
 
\begin{proposition} 
Let $G$ be a group satisfying \eqref{2} and the equivalent properties of Proposition 5.2.
We assume that $n \geq 2$.
\begin{enumerate}
\item For all $i$ in $\{2,\cdots,n\},$ there exists $\lambda_i \in \F_p-\{0\}$ such that $\ell_{i+1,i}=\lambda_i \, \ell_{2,1}$.\
Therefore, one can choose a basis of $D(G)$ in which the matrix $\Phi(y)$ reads as follows: 
$$ \Phi(y)= \begin{pmatrix}
 1 & 0 &  0 &\cdots& 0 \\
 \ell(y) & 1& 0 & \cdots  & 0 \\
 \ell_{3,1}(y)  & \ell(y) & \cdots&0 &0 \\
 \ell_{i,1}(y) & \ell_{i,2}(y) & \cdots & 1& 0\\
 \ell_{n,1}(y) & \ell_{n,2}(y) & \cdots & \ell(y) & 1
\end{pmatrix}$$
where $\ell$ is a nonzero linear form from $G/D(G)$ to $\F_p$.
\item Furthermore, $n \leq p$.
\end{enumerate}
\end{proposition}

\noindent \textbf{Proof:}
As $G/D(G)$ is abelian, for all $y$ and $y'$ in $V$, $\Phi(y)\Phi(y')=\Phi(y')\Phi(y)$. 
Then, for all $i$ in $\{ 1, \cdots, n-2\}$, the identification of the coefficients situated on the $(i+2)$-th line and the $i$-th column in the matrices $\Phi(y) \Phi(y')$ and $\Phi(y') \Phi(y)$ reads:
$$\ell_{i+2,i} (y)+\ell_{i+1,i}(y')\, \ell_{i+2,i+1}(y) +\ell_{i+2,i}(y')=
\ell_{i+2,i} (y')+\ell_{i+1,i}(y) \,\ell_{i+2,i+1}(y') +\ell_{i+2,i}(y)$$
Therefore, for all $y$ and $y'$ in $G/D(G)$, $\ell_{i+1,i}(y') \,\ell_{i+2,i+1}(y)=\ell_{i+1,i}(y') \,\ell_{i+2,i+1}(y)$.
As $\ell_{i+1,i}$ and $\ell_{i+2,i+1}$ are nonzero linear forms, it follows that $Ker\, \ell_{i+1,i}= Ker\, \ell_{i+2,i}$. Then, $\ell_{i+1,i}$ and $\ell_{i+2,i+1}$ are homothetic.
It implies that, for all $i$ in $\{2,\cdots,n\},$ there exists $\lambda_i \in \F_p-\{0\}$ such that $\ell_{i+1,i}=\lambda_i \, \ell_{2,1}$. We eventually replace the basis of $D(G)$: $(g_i)_{1\leq i \leq n}$ with 
$(\frac{1}{\lambda_{i}}\, g_i)_{1\leq i \leq n}$ and denote $\ell_{2,1}$ by $\ell$. In this new basis, the matrix $\Phi(y)$ reads as expected and the first point is proved.\\
\indent We now work with a basis of $D(G)$ in which the matrix $\Phi(y)$ reads as in the first point. We take some $y_0$ in $G/D(G)$ such that $\ell(y_0)\neq 0$.
One checks that $n$ is the smallest integer $m \geq 1$ such that $(\Phi(y_0)-I_n)^m=0$,
where $I_n$ denotes the identity matrix of size $n$. But, as $G/D(G)$ has exponent $p$, then $(\Phi(y_0)-I_n)^p= \Phi(y_0)^p-I_n=0$. It follows that $p \geq n$. 
 $\square$

 \begin{remarque}
In the situation exposed in Proposition 5.2, that is to say in the case where each linear form $\ell_{i+1,i}$ in $\Phi(y)$ is non identically zero, the representation $\phi$ is said to be indecomposable, i.e. if $D(G)=D(G)_1 \bigoplus D(G)_2$, where $D(G)_1$ and $D(G)_2$ are two $\F_p$-subvectors spaces of $D(G)$ stable by $\phi$, then the $D(G)_i$'s are trivial (left as an exercise to the reader).
Nevertheless, the converse is false, i.e. the representation $\phi$ can be indecomposable without the linear forms $\ell_{i+1,i}$'s being all nonzero.
\end{remarque}

\subsection{A group-theoretic characterization for big actions with $f_i \in \Sigma_{i+1}-\Sigma_i$.}
\indent In the sequel, we study the filtration defined by the $(\Lambda_i(G))_{i \geq 0}$ in the special case of a big action $(C,G)$ whose $G_2$ is $p$-elementary abelian. Since $G_2$ coincides with $D(G)$ (see \cite{MR}, Theorem 2.6), note that such a group $G$ systematically satisfies condition \eqref{2}. We now investigate the case where the group $G$ satisfies the equivalent properties of Proposition 5.2. In particular, we show that these group-theoretic conditions characterize the big actions with a $p$-elementary abelian $G_2$ and such that each $f_i$ lies in $\Sigma_{i+1}-\Sigma_i$.
The final section will be devoted to explicit families of big actions satisfying these properties.
Throughout this section, the notations concerning big actions are those fixed in section 3.2.

\begin{theoreme} 
Let $(C,G)$ be a big action with $G_2 \simeq (\Z / p\Z)^n$, for $n \geq 2$, and such that the group $G$ satisfies the equivalent properties of Proposition 5.2.
 Then, $p \geq n+1 \geq 3$. Furthermore, for all $i$ in $\{1,\cdots,n\}$, $m_i=1+ i \, p^{s_1}$. In particular, $f_i \in \Sigma_{i+1}-\Sigma_i.$
Moreover, $v=s_1+1$. In this case, $$\frac{|G|}{g}=\frac{2\,p}{p-1} \, 
\frac{p^n\, (p-1)^2}{n\,p^n\,(p-1)+1-p^n} >\frac{2\,p}{p-1}$$
\end{theoreme}

\noindent \textbf{Proof}: For a fixed $n$, we prove by induction on $i$ that for all $i$ in $\{1,\cdots,n\}$ such that $i \leq p-1$,  $m_i=1+ i \, p^{s_1}$.
By the way, we show that $n \leq p-1$. Indeed, we cannot propagate the induction when $i=p-1$ and $n=p$.

The first step of the induction derives from the definition of $m_1$. 
Then, we consider some integer $i$ such that $2 \leq i \leq n$ and $i \leq p-2$ and assume that the proposition is true for all $j \leq i-1$. As seen in section 2.4, we can write: 
\begin{equation} \label{3}
\forall \, y \in V,\quad  \Delta_y(f_i):=f_i(X+y)-f_i(X)= \sum_{j=1}^{i-1}\, \ell_{j,i}(y) \,f_j(X) \quad \mod \,  \wp(k[X]) 
\end{equation}
where the maps $\ell_{j,i}$ from $V$ to $\F_p$ refer to the coefficients of the matrix $L(y)$.
As the group $G$ satisfies the third condition of Proposition 5.2 which does not depend on the basis chosen for $D(G)$, it follows from Proposition 2.9 and Remark 2.10 that for all $i$ in $\{1,\cdots,n-1\}$, each $\ell_{i,i+1}$ is a nonzero linear form from $V$ to $\F_p$.

\begin{enumerate}
\item
\textit{We first prove that the function $f_i$ does not belong to $\Sigma_i$. }\\
Assume that $f_i$ lies in $\Sigma_i$ and apply Lemma 3.12 to $f(X):=\Delta_y(f_i)-\sum_{j=1}^{i-1}\, \ell_{j,i}(y) \,f_j(X)$ and $a_0:=m_{i-1}$.
By induction hypothesis, $m_{i-1}=1+(i-1)\,p^{s_1} \in \N-p\N$. Note that $X^{a_0}=X^{1+(i-1)\,p^{s_1}}$ lies in $\Sigma_i-
\Sigma_{i-1}$. We gather from Lemma 3.7.6 that no $X^{a_0p^r}$, with $r \geq 0$, belongs to $\Sigma_{i-1}$, so none of them can be found in $\Delta_y(f_i)$ which belongs to $\Sigma_{i-1}$, as $f_i$ lies in $\Sigma_i$ (cf. Lemma 3.9). 
 Besides, the property c imposed on $m_i$ by Definition 2.3.4 implies that, for any $y$ in $V$ such that
 $\ell_{i-1,i}(y) \neq 0$, $a_0=m_{i-1}$ is the degree of $\sum_{j=1}^{i-1}\, \ell_{j,i}(y) \,f_j(X)$.
 Such an element $y$ exists since $\ell_{i-1,i}$ is supposed to be a nonzero linear form. It follows that,
 when keeping the notation of Lemma 3.12, $f_{a_0}(X)=c_{m_{i-1}}(f_{i-1})\, \ell_{i-1,i}(y) X^{a_0}$, where
 $c_{m_{i-1}}(f_{i-1})\neq 0$ denotes
the coefficient of $X^{m_{i-1}}$ in $f_{i-1}$. As $p$ does not divide $a_0$, we gather from Lemma 3.12 that $f_{a_0}(X)$ is identically zero, which contradicts $\ell_{i-1,i}(y) \neq 0$. Therefore 
$f_i$ does not belong to $\Sigma_i$. In particular, as $\Sigma_2 \subset \Sigma_i$, $f_i$ does not belong to $\Sigma_2$.
Accordingly, we can define an integer $a \leq m_i$ such that $X^a$ is the monomial of $f_i$ with highest degree
which does not lie in $\Sigma_2$. Since $f_i$ is assumed to be reduced mod $\wp(k[X])$, 
$a \not \equiv 0$ mod $p$.
  
\item \textit{We now prove that $a-1 \geq 1+(i-1)\, p^{s_1}$.}\\
Assume that $a-1< 1+(i-1)\, p^{s_1}$ and apply Lemma 3.12 to $f(X):=\Delta_y(f_i)-\sum_{j=1}^{i-1}\, \ell_{j,i}(y) \,f_j(X)$ and $a_0:=m_{i-1}=1+(i-1)\,p^{s_1} \in \N-p\N$. The proof works as above except that we now have to determine the monomials of $f_i$ which could produce some $p$-powers of $X^{a_0}$ in $\Delta_y(f_i)$. As $a-1<a_0$, they must be searched for among the monomials of $f_i$ with degree strictly greater than $a$. But, by definition of $a$, such monomials belongs to $\Sigma_2$, so give monomials in $\Delta_y(f_i)$ which are in $\Sigma_1$, whereas $X^{a_0}=X^{1+(i-1)\,p^{s_1}}$ lies in $\Sigma_i-\Sigma_{i-1}$, with $i \geq 2$. Just as in the first point, we can conclude that, for any $y$ in $V$ such that $\ell_{i-1,i}(y) \neq 0$,  $f_{a_0}(X)=c_{m_{i-1}}(f_{i-1})\, \ell_{i-1,i}(y) X^{a_0}$, which leads to the same contradiction as above.

\item \textit{We show that $p$ divides $a-1$.}\\
Assume that $p$ does not divide $a-1$. We first suppose that $a-1> 1+(i-1)\, p^{s_1}$ and apply Lemma 3.12
to $f(X):=\Delta_y(f_i)-\sum_{j=1}^{i-1}\, \ell_{j,i}(y) \,f_j(X)$ and $a_0:=a-1 \in \N-p\N$. As explained above, the monomials in $f_i$ with degree stricly greater than $a$, produce in $\Delta_y(f_i)$ monomials which are in $\Sigma_1$. But, as $p$ does not divide $a-1$, the monomial $X^{a-1}$ cannot belong to $\Sigma_1$: otherwise, $a-1=1$, which contradicts $a-1> 1+(i-1)\, p^{s_1}$, with $i \geq 2$. So the only $p$-power of $X^{a-1}$ that occur in $\Delta_y(f_i)$ comes from the monomial $X^a$ of $f_i$: it is $c_{a}(f_{i})\,a \, y\,  X^{a-1}$, where $c_{a}(f_{i})\neq 0$ denotes the coefficient of $X^{a}$ in $f_{i}$.
Besides, $X^{a-1}$ does not occur in $\sum_{j=1}^{i-1}\, \ell_{j,i}(y) \,f_j(X)$ whose degree is at most
$1+(i-1)\, p^{s_1}< a-1$.
We gather from Lemma 3.12 that $f_{a_0}(X)=c_{a}(f_{i})\,a \, y\,  X^{a_0}$ is 
identically zero. It implies that $V=\{0\}$, which is excluded for a big action.
Accordingly, $a-1= 1+(i-1)\,p^{s_1}$. The equality of the leading coefficients in \eqref{3} implies that for all $y$ in $V$, $ \ell_{i-1,i}(y)=\frac{a\, c_a(f_i)}{c_{m_{i-1}}(f_{i-1})}\,y$. So the kernel of the linear form
$\ell_{i-1,i}$ is reduced to $\{0\}$ and $v \leq 1$, which contradicts Lemma 3.11.
\indent Accordingly, $p$ divides $a-1$. Thus, we can write $a:= 1+\lambda \, p^t$, with $t >0$, $\lambda$ prime to $p$ and $\lambda \geq 2$ because of the definition of $a$.

\item \textit{Put $j_0:= a-p^t=1+(\lambda-1)\, p^t$. We prove that $j_0 = 1+(i-1)\, p^{s_1}$.} \\
Indeed, if  $j_0<  1+(i-1) \, p^{s_1}$, then $a=j_0+p^t < 1+(i-1)\,p^{s_1}+p^t.$ 
Using the second point, we get: $1+(i-1)\, p^{s_1} < a= 1+\lambda \, p^t < 1+p^t+(i-1)\, p^{s_1}.$ 
If $s_1 - t \geq 0 $, it implies $(i-1)\, p^{s_1-t} < \lambda < 1+(i-1) \,p^{s_1-t}$ with $p^{s_1-t} \in \N$, which is impossible. So, $s_1-t \leq -1$. In this case, as $i-1 < p$, the inequality  $1+ \lambda \, p^{t} < 1 +(i-1)\, p^{s_1}+p^t$ yields: $\lambda-1 < (i-1)\, p^{s_1-t} < p^{1+s_1-t} \leq 1$,  which contradicts $\lambda \geq 2$. As a consequence, $j_0 \geq 1+(i-1)\, p^{s_1}$.\\
  \indent We now prove that $j_0 = 1+(i-1)\, p^{s_1}.$ 
  Assume that $j_0 > 1+(i-1)\, p^{s_1}$ and apply Lemma 3.12 to $f(X):=\Delta_y(f_i)-\sum_{j=1}^{i-1}\,
  \ell_{j,i}(y) \,f_j(X)$ and $a_0:=j_0 \in \N-p\N$. No $p$-power of $X^{j_0}$ can be found in
  $\sum_{j=1}^{i-1}\, \ell_{j,i}(y) \,f_j(X)$ whose degree is at most  $1+(i-1)\, p^{s_1} <j_0$.
  It follows that the monomials $X^{j_0p^r}$ have to be searched for in $\Delta_y(f_i)$. Then,
  the same argument as in the proof of Theorem 3.13 allows to write:
  $f_{a_0}(X)= T(y)\, X^{j_0}$ with $T(y):= \sum _{b=j_0+1}^a  \, c_b(f_i) \,\binom b{j_0} \,y^{b-j_0}$,
  where $c_b(f_i)$ denotes the coefficient of $X^b$ in $f_i(X)$. This entails the same contradiction with
  Lemma 3.11 as in the proof of Theorem 3.13. Therefore, $j_0 = 1+(i-1)\, p^{s_1}$.

\item \textit{We gather that $v=t+1$.} \\
Indeed, since $j_0= 1+(i-1)\, p^{s_1}= deg\, f_{i-1}$ , the equality of the corresponding coefficients in \eqref{3} reads: $T(y)=\ell_{i-1,i}(y) \, c_{m_{i-1}} (f_{i-1})$, which holds for all $y$ in $V$. Put $\tilde{T}:=\frac{T}{c_{m_{i-1}} (f_{i-1})}$. It has the same degree as $T$ and satisfies $\tilde{T}(y)=\ell_{i-1,i}(y) \in \F_p$, for all $y$ in $V$. It follows that $\tilde{T}^p -\tilde{T}$ is identically zero on $V$, so $v \leq t+1$. Using the same argument as in the proof of Theorem 3.13, we prove that $v \leq t$ contradicts Lemma 3.11. We gather that $v=t+1$.

\item  \textit{We prove that $s_1=t$. It follows that $v=s_1+1$ and $a=1+i\, p^{s_1}$, which requires $p>n \geq 2$.}\\
 \indent As $j_0=1+(i-1)\, p^{s_1}$, then $a=j_0+p^t=1+(i-1)\, p^{s_1}+p^t$.
  But, $a=1+\lambda \, p^t \geq 1+2 \, p^t$.  From $i-1\leq p$, we gather that  $p^t \leq (i-1) \, p^{s_1} <
  p^{s_1+1}$.
  Therefore, $t \leq s_1$. To prove the equality, we focus on the big action $(C/H_i, G/H_i)$ as defined
  in Proposition 2.11. Since $v =t+1$, then $|G/H_i|=p^{i+v}= p^{i+t+1}$. Besides, as $m_i \geq a=j_0+p^t=1+(i-1)
  \, p^{s_1}+ p^t \geq 1+p^{s_1}+p^t$, $$g_{C/H_i} \geq  \frac{(p-1)}{2}\, p^{i-1}\, (m_{i}-1)\geq
  \frac{(p-1)}{2}\, (p^{i-1+s_1}+p^{i-1+t})$$
   If $u:=s_1-t \geq 1$, the lower bound for the genus becomes: $$g_{C/H_i} \geq \frac{(p-1)}{2}\,
   p^{t+i-1}\, 
(p^u+1) \geq \frac{(p-1)}{2}\, p^{t+i-1}\,(p+1)$$ This contradicts condition $(N)$ insofar as: $$\frac{|G/H_i|}{g_{C/H_i}} \leq \frac{2\: p}{p-1} \, \frac{p^{i+t}}{p^{i-1+t} \, (p+1)} =\frac{2\: p}{p-1} \, \frac{p}{p+1} <\frac{2\: p}{p-1}$$ 
Therefore, $s_1=t$, so $v=s_1+1$ and $a=1+(i-1)\, p^{s_1}+p^t= 1+i\, p^{s_1}.$
Note that we find: $\lambda=i$. As $\lambda$ is supposed to be prime to $p$ and as $2 \leq i \leq n$ and $i \leq p-1$, it requires that $p >n \geq 2$.

\item \textit{We conclude that $m_i=a=1+i\, p^{s_1}$.}\\
Assume $a < m_i$. Then, by definition of $a$, there exists an integer $r \geq 0$ such that $m_i=1+p^r$. Thus, we get: $m_i=1+p^r> a=1+i\, p^{s_1} \geq 1+2\, p^{s_1}$. As $p \geq 3$, this implies $r \geq s_1+1$.
We gather a new lower bound for the genus of $C/H_i$, namely: 
$$g_{C/H_i} \geq \frac{(p-1)}{2}\,(p^{s_1} +p^{i-1}\, (m_i-1))=
\frac{(p-1)}{2}\,(p^{s_1} +p^{i-1+r})\geq \frac{(p-1)}{2}\,(1+p^{i+s_1})$$
As $|G/H_i|=p^{i+s_1+1}$, it follows that $\frac{|G/H_i|}{g_{C/H_i}} \leq \frac{2\, p}{(p-1)}\, \frac{p^{i+s_1}}{1+p^{i+s_1}} < \frac{2\, p}{(p-1)} $, which contradicts condition $(N)$ for the big action $(C/H_i,G/H_i)$. 
\indent Accordingly, $m_i=a=1+i\, p^{s_1}$, which completes the induction.
\end{enumerate}
To conclude,  we compute the genus by means of Corollary 2.7, namely
$g= \frac{p-1}{2} \, p^{s_1} \, \frac{n\,p^n \,(p-1)+1-p^n}{(p-1)^2}$.
It follows that $$\frac{|G|}{g}=\frac{2\,p}{p-1} \, 
\frac{p^n\, (p-1)^2}{n\,p^n\,(p-1)+1-p^n} \geq 
\frac{2\,p}{p-1} \, 
\frac{p^n\, (p-1)^2}{p^n\,(p-1)^2+1-p^n} > \frac{2\,p}{p-1} \qquad \quad \square$$

\begin{corollaire}
Let $(C,G)$ be a big action as in Theorem 5.6.
Let $G_{\infty,1}$ be the wild inertia subgroup of $Aut_k(C)$ at $\infty$.
Then, $G$ is equal to $G_{\infty,1}$.
\end{corollaire}

\noindent \textbf{Proof:}
As before, we denote by $L$ the function field of $C$ and by $k(X)$ the subfield of $L$ fixed by $D(G)$.
By \cite{MR} (Corollary 2.10), the pair $(C,G_{\infty,1})$ is a big action such that $D(G_{\infty,1})
=D(G)$. It follows that $G/D(G)$ is included in  $G_{\infty,1}/D(G_{\infty,1})$, both of them acting
 as a group of translations of $Spec\, k[X]$. In the same way as we define the
 representation $\phi: G/D(G) \rightarrow Aut(D(G))$ in section 2.2 (or more generally in section 6.1), consider a representation $\tilde{\phi}$ from $G_{\infty,1}/D(G_{\infty,1})$ to $Aut(G_{\infty,1})$.
Fix an adapted basis of $A$ and then, by duality, a basis of $D(G)$ in which, for all $y$ in $G/D(G)$ the matrix $\Phi(y)$ of the automorphism $\phi(y)$ is lower triangular. For all $y$ in $G_{\infty,1}/D(G_{\infty,1})$, call $\tilde{\Phi}(y)$ the matrix of the automorphism $\tilde{\phi}(y)$ in the same adapted basis. When restricted to $G/D(G)$, the two matrices coincides, i.e. if $y$ lies in $G/D(G) \subset G_{\infty,1}/D(G_{\infty,1})$, $\Phi(y)=\tilde{\Phi}(y)$. As a consequence, the group $G_{\infty,1}$ also satisfies the third condition 
of Proposition 5.2. Therefore, by Theorem 5.6,  $\frac{|G_{\infty,1}|}{|D(G_{\infty,1})|}=s_1+1=\frac{|G|}{|D(G)|}$. So, $G=G_{\infty,1}$.
$\square$
\medskip

We conclude this section by showing that the big actions $(C,G)$ such that $G$ satisfies the equivalent properties of Proposition 5.2 are exactly those with $f_i \in \Sigma_{i+1}-\Sigma_i$.

\begin{theoreme}
Let $(C,G)$ be a big action with $G_2 \simeq (\Z/p\Z)^n$, for $n \geq 2$. 
We keep the notation defined in sections 4.1 and 4.2.
Then, the following assertions are equivalent.
\begin{enumerate}
\item For all $i$ in $\{1,\cdots,n\}$, the function $f_i$ lies in $\Sigma_{i+1}$-$\Sigma_i$. 
\item  The group $G$ satisfies the equivalent properties of Proposition 5.2.
\end{enumerate}
In this case, $n \leq p-1$, $m_i=1+i\, p^{s_1}$ for all $i$ in $\{1,\cdots,n\}$, $v=s_1+1$ and $\frac{|G|}{g}=\frac{2\,p}{p-1} \, 
\frac{p^n\, (p-1)^2}{n\,p^n\,(p-1)+1-p^n}$. Moreover, the upper ramification groups of $G_2$ coincide with the subgroups $\Lambda_i(G)$ studied in section 5.1. More precisely, following the notation of section 2.3.2, $(G_2)^{\nu_{i}}=\Lambda_{n-i}(G)$ for all $i$ in $\{0, \cdots,n\}$.
\end{theoreme}

\noindent \textbf{Proof:} 
The implication from $2$ to $1$ comes from Theorem 5.6 which also shows that, in this case, $n \leq p-1$, $m_i=1+i\, p^{s_1}$, for all $i$ in $\{1,\cdots,n\}$, $v=s_1+1$ and $\frac{|G|}{g}=\frac{2\,p}{p-1} \, 
\frac{p^n\, (p-1)^2}{n\,p^n\,(p-1)+1-p^n}> \frac{2\,p}{p-1}$. Conversely, assume that the second assertion is
satisfied. We prove by induction on $i$ that, for all $i$ in $\{1,\cdots,n-1\}$, the linear form $\ell_{i,i+1}$ is nonzero. Then, by Proposition 2.9, Remark 2.10 and Remark 5.3, we gather that the group $G$ satifies the third condition of Proposition 5.2.
We first study the case $i=1$ and consider the big action $(C/H_2,G/H_2)$, as defined in Proposition 2.11, i.e. the big action whose curve $C/H_2$ is parametrized by $W_j^p-W_j=f_j(X)$, with $1 \leq j \leq 2$. By hypothesis,
$f_2$ does not lie in $\Sigma_2$. We infer from Proposition 2.13 that the representation $\rho$ associated with $(C/H_2,G/H_2)$ is non trivial. Then, the linear form $\ell_{1,2}$ is nonzero.
We now take $i \geq 2$ and assume that the property is true for all $j \leq i$. It means that, for all $i$ in $\{1,\cdots,i-1\}$, the linear form $\ell_{j,j+1}$ is nonzero. Then, by Theorem 5.6, for all $j$ in $\{1,\cdots,i-1\}$, $m_j=1+j \, p^{s_1}$ and 
$v=s_1+1$. We now write condition $(N)$ for the big action $(C/H_{i+1}, G/H_{i+1})$ as defined in Proposition 2.11, that is to say the big action parametrized by $W_j^p-W_j=f_j(X)$, with $1 \leq j \leq i+1$. As $|G/H_{i+1}|=p^{v+i+1}=p^{s_1+i+2}$ and $g_{C_{H_{i+1}}}=\frac{p-1}{2} \, \{(\sum_{j=1}^{i}\, j \, p^{s_1+j-1}) +p^i \, (m_{i+1}-1)\} $, we gather that the inequality $\frac{|G_{H_{i+1}}|}{g_{C_{H_{i+1}}}}>\frac{2\,p}{p-1}$ is equivalent to the following condition on $m_{i+1}$:
\begin{equation} \label{4} 
m_{i+1} <p^{s_1+1}-(\sum_{j=1}^i \, j\,  p^{s_1+j-1-i}) +1=p^{s_1} (p-1)+\sum_{j=2}^i (p-(j+1))\, p^{s_1+j-i-1} +(p-1) \, p^{s_1-i} +1 
\end{equation}
We now assume that $\ell_{i,i+1}$ is the null linear form. Then, for all $y$ in $V$,
$\Delta_y(f_i)=\sum_{j=1}^{i-1}\, \ell_{j,i+1} (y) \, f_j(X)$ mod $\wp(k[X])$. This ensures that the function field of the curve $\mathcal{C}: W_j^p-W_j=f_j(X)$, with $1 \leq j \leq i+1$ and $j \neq i$, is a Galois extension of $k(X)$ whose group $\mathcal{H}$ is isomorphic to $(\Z/p\Z)^i$ and, as usual, the group of translations by $V$ extends to an automorphism group of $\mathcal{C}$, say $\mathcal{G}$, with the following exact sequence:
$$0 \longrightarrow \mathcal{H} \longrightarrow \mathcal{G} \longrightarrow V \longrightarrow 0$$
We compute the quotient $\frac{|\mathcal{G}|}{g_{\mathcal{C}}}$. As $|\mathcal{G}|=p^{s_1+i+1}$ and 
$g_{\mathcal{C}}= \frac{p-1}{2} \, \{(\sum_{j=1}^{i-1}  j\, p^{s_1+j-1}) +\, p^{i-1} \, (m_{i+1}-1)\}$, one can check that 
$\frac{|\mathcal{G}|}{g_{\mathcal{C}}} > \frac{2\,p}{p-1}$ if and only if 
\begin{equation} \label{5}
m_{i+1} <p^{s_1+1} -\sum_{j=1}^{i-1} \, j \, p^{s_1+j-i} +1=
p^{s_1}\, (p-1) + \sum_{j=1}^{i-2} \, (p-(j+1)) \, p^{s_1+j-i} +(p-1)\, p^{s_1-i+1}+1 
\end{equation}
The condition \eqref{5} is verified since it is implied by \eqref{4}.
It follows that $(\mathcal{C}, \mathcal{G})$ is a big action. By Theorem 3.13, it implies that the $i$-th function: $f_{i+1}$,
lies in $\Sigma_{i+1}$, which contradicts the hypothesis $f_{i+1} \in \Sigma_{i+2}-\Sigma_{i+1}$.
Therefore, $\ell_{i,i+1}$ is a nonzero linear form, which completes the induction and prove the equivalence between $1$ and $2$.\\
\indent We now prove the last statement on the upper ramification filtration of $G_2$. Starting from a given 
adapted basis of $A$: $\{\overline{f_1(X)}, \cdots, \overline{f_n(X)}\}$, we get, by duality with respect to the Artin-Schreier pairing, a basis of $G_2$, say $\{g_1,\cdots,g_n\}$. As proved, for all $i$ in $\{1, \cdots,n\}$, $m_i=1+i\, p^{s_1}$, the jumps in the upper ramification of $A$, as defined in section 2.3.1, are:
$\mu_i=m_{i+1}=1+(i+1)p^{s_1}$, for all $i$ in $\{0,\cdots,n-1\}$. Put $\mu_n:=1+m_n$. Then,  $A^{\mu_0}=\{0\}$ and, for all $i$ in $\{1,\cdots,n\}$, $A^{\mu_i}$ is the $\F_p$-subvector space of $A$ generated by $\overline{f_1(X)}, \cdots, \overline{f_i(X)}$. By duality (see Proposition 2.5), $(G_2)^{\nu_n}=(G_2)^{\mu_n}=\{e\}=\Lambda_0(G)$ and, for all $i$ in $\{0,\cdots,n-1\}$, 
$(G_2)^{\nu_i}=(G_2)^{\mu_i}$ is the $\F_p$-subvector space of $G_2$ generated by $g_{i+1}, \cdots, g_n$,
which is precisely $\Lambda_{n-i}(G)$, as seen in Proposition 5.2. $\square$

\section{Examples.}
 \indent We conclude this paper with some examples illustrating the special case of big actions described in Theorem 5.8, namely the big actions $(C,G)$ with a $p$-elementary abelian $G_2$ such that each $f_i$ lies in $\Sigma_{i+1}-\Sigma_i$. Note that Theorem 5.8 is twofold: on the one hand, it gives a group-theoretic
characterization of $G$ (cf. 5.8.2) and, on the other hand, it displays a dual point of view related to the parametrization of the cover (cf. 5.8.1). When studying the special family explicitely constructed via equations in Proposition 6.1, the second point of view naturally dominates in the proof. On the contrary, when exploring a universal family as in section 6.2, we are lead to combine both aspects.\\

\textit{Notation.} The notations concerning big actions are those fixed in section 3.2. Moreover, let $W(k)$ be the ring of Witt vectors with coefficients in $k$. Then, for any $\sigma \in k$, we denote by $\tilde{\sigma}$ the Witt vector $\tilde{\sigma}:=(\sigma,0,0,\cdots) \in W(k)$. For any $S(X):=\sum_{i=0}^{s} \, \sigma_i \, X^i \in k[X]$, we denote by $\tilde{S}(X)$ the polynomial $\sum_{i=0}^s \, \tilde{\sigma_i} \, X^i \in W(k)[X]$. 

\subsection{A special family.}
\subsubsection{Case $s_1=1$.}
\indent Let $p \geq 3$ and $1\leq n \leq p-1$. We first exhibit a special family of big actions $(C,G)$ which satisfy the conditions of Theorem 5.8 with $s_1=1$ and so, $v=dim_{\F_p} V=2$. We shall distinguish the cases $n <p-1$ and $n=p-1$. When investigating the properties of the corresponding group $G$, we show, among others, that $G$ is a capable group (see Definition 6.8) as studied by \cite{Ha} and \cite{BT}. 

\begin{proposition}
Let $p \geq 3$.
Let $S(X):=\wp(X)$ and $Q(X):=\wp(S(X))$. Call $V$ the $\F_p$-vector space $V$ consisting of the set of zeroes of the
polynomial $Q$. Then, $V \simeq (\Z/p \Z)^2$.
\begin{enumerate}
\item Let $n$ in $\{1,\cdots, p-2\}$. For all $i$ in $\{1,\cdots,n\}$, we denote
$$g_i(X):= \frac{S(X)^{i+1}}{(i+1)!}$$
Let $f_i:=red(g_i)$ be the reduced representative of $g_i$, as defined in the introduction.
Let $C[n]$ be the curve parametrized by the $n$ Artin-Schreier equations:
 $W_i^p-W_i=f_i(X)$, for $ 1 \leq i \leq n$.
 Let $K_n:=k(C[n])$ be the function field of $C[n]$ and $H[n] \simeq (\Z/p\Z)^n$ be the Galois group of $K_n/k(X)$. Then,
 the group of translations of the affine line: $\{\tau_y:\, X \rightarrow X+y , \, y \in V\}$ extends to an automorphism
 $p$-group of $C[n]$, say $G[n]$, such that we get the exact sequence:
$$ 0 \longrightarrow H[n] \simeq (\Z/p \,\Z)^{n} \longrightarrow G[n] \longrightarrow V\simeq (\Z/p \Z)^2
  \longrightarrow 0$$
In this case, the pair $(C[n],G[n])$ is a big action with $G[n]_2 \simeq (\Z/p \,\Z)^{n}$. Moreover, this big action satisfies the conditions of Theorem 5.8 with $s_1=1$.
\item We now study the case $n=p-1$.
We define $g_{p-1}(X) \in k[X]$ as the reduction
mod $p$ of the polynomial $$\frac{1}{p!} \,((X^p-X)^p-X^{p^2}+X^p) \in \, W(k)[X]$$
 Let $f_{p-1}$ be the reduced representative of $g_{p-1}$. 
 Let $C[p-1]$ be the curve parametrized by the $p-1$ Artin-Schreier equations:
 $W_i^p-W_i=f_i(X)$, for $ 1 \leq i \leq p-1$, where the $p-2$ first $f_i$'s are defined as in the
first case.
 Let $K_{p-1}:=k(C[p-1])$ be the function field of $C[p-1]$ and $H[p-1] \simeq (\Z/p\Z)^{p-1}$ be the Galois group of
$K_{p-1}/k(X)$. Then, the group of translations of the affine line: $\{\tau_y:\, X \rightarrow X+y, \, y \in V\}$ extends to an automorphism
 $p$-group of $C[p-1]$, say $G[p-1]$, with the following exact sequence:
$$ 0 \longrightarrow H[p-1] \simeq (\Z/p \,\Z)^{p-1} \longrightarrow G[p-1] \longrightarrow V\simeq (\Z/p \Z)^2
  \longrightarrow 0$$
In this case, the pair $(C[p-1],G[p-1])$ is a big action with $G[p-1]_2 \simeq (\Z/p \,\Z)^{p-1}$. Moreover, this big action satisfies the conditions of Theorem 5.8 with $s_1=1$.

\end{enumerate}
\end{proposition}

\noindent \textbf{Proof:} 
Using Proposition 2.11, we first observe that the second case implies the first one, when excluding the last equation:
$W^p_{p-1}-W_{p-1}=f_{p-1}(X)$. Therefore, it is sufficient to prove the second point.\\
\indent Fix $y \in V$. We begin by calculating $\Delta_y(g_i)$ for $1 \leq i \leq p-2$. So, take $i$ in $\{1, \cdots, p-2\}$. One first shows that 
 $$  \Delta_y(g_i)=g_i(X+y)-g_i(X)= \sum_{j=1}^{i-1} \, \frac{S(y)^{i-j}}{(i-j)!} \,
 g_j(X) +g_i(y)+ \frac{S(y)^i}{i!} \, S(X)$$ 
where the first sum is empty when $i=1$. Since $S(y)$ lies in $\F_p$ for all $y$ in $Z(Q)=V$, one gets:
 $$\Delta_y(g_i)=\sum_{j=1}^{i-1} \, \frac{S(y)^{i-j}}{(i-j)!} \,
 g_j(X)+g_i(y) + \wp(\frac{S(y)^i}{i!} \, X).$$
As $k$ is an algebraically closed field, $g_i(y)=0$ mod $\wp(k[X])$. We gather that
  $\Delta_y(g_1)=0$ mod $\wp(k[X])$ and that, for all $i$ in $\{2,\cdots, p-2\}$,
 $\Delta_y(g_i)= \sum_{j=1}^{i-1} \, \ell_{j,i}(y) \,g_j(X)$ mod $\wp(k[X])$ with
 $\ell_{j,i}(y):= \frac{S(y)^{i-j}}{(i-j)!} \in \F_p$.
 Since $g_i(X)-f_i(X)$ lies in $\wp(k[X])$ and since each $\ell_{j,i}(y)$ belongs to $\F_p$, it follows that $$\forall \, i \in \{1, \cdots, p-2\},  \quad \Delta_y(f_i)= \sum_{j=1}^{i-1} \, \ell_{j,i}(y) \,f_j(X) \quad \mod \, \wp(k[X]) \quad \mbox{with} \quad \ell_{j,i}(y):= \frac{S(y)^{i-j}}{(i-j)!}$$
\indent Now fix $y$ in $V$ and calculate $\Delta_y(g_{p-1})$. 
As $X^p-X=S(X)$ mod $p$, we first notice that $(X^p-X)^p=\tilde{S}(X)^p$ mod $p^2  W(k)[X]$. It follows that $g_{p-1}$ can also be seen as the reduction mod $p$ of the polynomial: $\frac{1}{p!} \,(\tilde{S}(X)^p-X^{p^2}+X^p) \in W(k)[X]$. 
By the same token, from $S(X+y)=S(X)+S(y)$ mod $p$, we gather that 
$\tilde{S}(X+\tilde{y})^p=(\tilde{S}(X)+\tilde{S}(\tilde{y}))^p$ mod $p^2  W(k)[X]$.
 It follows that $$\tilde{S}(X+\tilde{y})^p-\tilde{S}(X)^p-\tilde{S}(\tilde{y})^p=
 \sum_{i=1}^{p-1} \, \binom pi \, \tilde{S}(X)^i \, \tilde{S}(\tilde{y})^{p-i} \, \mod \, p^2 W(k)[X]$$
Likewise,
$$(X+\tilde{y})^p-X^p-\tilde{y}^p-(X+\tilde{y})^{p^2}
+X^{p^2}+\tilde{y}^{p^2}=\sum_{i=1}^{p-1} \, \binom pi \, (X^i\, \tilde{y}^{p-i}-
X^{pi}\, \tilde{y}^{p(p-i)}) \, \mod    \, p^2\, W(k)[X]$$
  Then, we obtain the following equalities :
  $$\frac{1}{p!}
(\tilde{S}(X+\tilde{y})^p-\tilde{S}(X)^p-\tilde{S}(\tilde{y})^p+(X+\tilde{y})^p-X^p-\tilde{y}^p-(X+\tilde{y})^{p^2}
+X^{p^2}+\tilde{y}^{p^2})$$
$$= \sum_{i=1}^{p-2} \,
\frac{\tilde{S}(\tilde{y})^{p-i-1}}{(p-i-1)!} \, \frac{\tilde{S}(X)^{i+1}}{(i+1)!} \, + \frac{\tilde{S}(\tilde{y})^{p-1}}{(p-1)!} \, \tilde{S}(X)\,+\sum_{i=1}^{p-1} \, \frac{\binom pi}{p!} \, (X^i\, \tilde{y}^{p-i}-
X^{pi}\, \tilde{y}^{p(p-i)}) \, \mod \, p W(k)[X]$$
$$=\sum_{i=1}^{p-2} \, \frac{S(y)^{p-i-1}}{(p-i-1)!} \, \frac{S(X)^{i+1}}{(i+1)!} \,+
\frac{S(y)^{p-1}}{(p-1)!}\, S(X)+\sum_{i=1}^{p-1} \, \frac{(-1)^i}{i} \, (X^i\, y^{p-i}-X^{ip}\, y^{p(p-i)})
\, \mod \, p W(k)[X]$$
since the kernel of the map:$ 
\left\{
\begin{array}{lc}
W(k) \rightarrow k \\
(a_0,a_1,\cdots) \rightarrow a_0
\end{array}
\right.
$ 
is $p W(k)$.
From $S(y) \in \F_p$ , we infer:
$$\Delta_y(g_{p-1})=\sum_{i=1}^{p-2} \, \frac{S(y)^{p-i-1}}{(p-i-1)!} \,g_i(X) +\wp(\frac{S(y)^{p-1}}{(p-1)!}\, X)+\wp(\sum_{i=1}^{p-1} \, \frac{(-1)^{i+1}}{i} \, X^i\, y^{p-i})+g_ {p-1}(y)$$
It follows that  $\Delta_y(g_{p-1})= \sum_{i=1}^{p-2} \, \ell_{i,p-1}(y)\, g_i(X)$ mod $\wp(k[X])$, with 
$\ell_{i,p-1}(y)= \frac{S(y)^{p-1-i}}{(p-1-i)!} \in \F_p$. 
Since $g_i-f_i \in \wp(k[X])$ and
$\ell_{i,p-1}(y) \in \F_p$, $$\Delta_y(f_{p-1})= \sum_{i=1}^{p-2} \, \ell_{i,p-1}(y)\, f_i(X)\quad \mod \, \wp(k[X])\quad \mbox{with} \quad \ell_{i,p-1}(y)= \frac{S(y)^{p-1-i}}{(p-1-i)!}$$  
By Galois Theory, this ensures that the group $G[p-1]$ is well-defined.
Furthermore, it is easy to check that for all $i$ in $\{1,\cdots,p-1\}$, $deg\, f_i=1+i\,p$.
In this case, the same computation as in the end of the proof of Theorem 5.6 shows
 that $\frac{|G[p-1]|}{g_{C[p-1]}}=\frac{2\,p}{p-1} \, 
\frac{p^{p-1}\, (p-1)^2}{(p-1)\,p^{p-1}\,(p-1)+1-p^{p-1}}$, which proves that the pair $(C[p-1],G[p-1])$ is a big action.
To conclude, note that for all $i$ in $\{1, \cdots,p-2\}$ and for all $y$ in $V$, 
 $\ell_{i,i+1}(y)=S(y)$, which proves that $\ell_{i,i+1}$ is a nonzero linear form from $V$ to $\F_p$. 
 Therefore, because of Remarks 2.10 and 6.3, $G[p-1]$ satisfies the third assertion of Proposition 5.2 and then the conditions of Theorem 5.8. $\square$

\begin{remarque}
The preceding proof shows that, in the case of Proposition 6.1, 
$$\forall \, y \in V, \, \, \forall \, i \in \{2,\cdots,p-1\} \quad \mbox{and} \quad \forall \, j \in \{1,\cdots, i-1\}, \quad \ell_{j,i}(y)=\frac{S(y)^{i-j}}{(i-j)!}$$
It follows that the matrix $L(y)$ defined in section 2.4 reads:
$L(y)= exp (S(y)\, J)=\sum_{i=0}^{n-1} \, \frac{(S(y)J)^i}{i!}$ 
where $J$ is the $n\times n$ nilpotent matrix:  $$J:= \begin{pmatrix}
 0 &1& 0 &\cdots& 0 \\
 0 & 0 & 1 & \cdots& 0 \\
 0 & 0 &0 &\cdots&0 \\
 0 & 0 &0 &0 & 1& \\
 0 & 0 & 0 & 0 & 0
\end{pmatrix}$$

\end{remarque}

\begin{proposition}
Let $(C[n],G[n])$ be the big action described in the first point of Proposition 6.1, i.e. with 
$n <p-1$. The notations are those introduced in Proposition 6.1.
\begin{enumerate}
\item Let $\sigma$ in $G[n]$. Let $y$ in $V$ such that $y:=\sigma(X)-X$.
Then,
$$\sigma [W]= ^tL(y) [W]+ X [\mathcal{R}(y)]+ [Z(y)]$$
where $^tL(y)$ denotes the transpose matrix of the upper triangular matrix $L(y)$ defined in section 2.4., $[W]:=^t[W_1,\cdots,W_n]$,  $[\mathcal{R}(y)]:=^t[\frac{S(y)}{1\,!}, \cdots, \frac{S(y)^n}{n\, !}]$, and $[Z(y)]:=^t[Z_1(y),\cdots,Z_n(y)]$, where, for all $i$ in $\{1,\cdots,n\}$,  $Z_i(y)$ is an element of $k$ which satisfies $\wp(Z_i(y))=g_i(y)$.
\item The group $G[n]$ has exponent $p$. 
\end{enumerate}
\end{proposition}

\noindent \textbf{Proof:}
\begin{enumerate}
\item For the need of the proof, it is more convenient to work with the non-reduced functions, namely the functions $g_i$'s.
However, we still write the equations: $W_i^p-W_i=g_i(X)$, without changing the notation of $W_i$.
As seen in the proof of Proposition 6.1, 
$$\forall \, y \in V, \quad  \forall \, i \in \{1,\cdots,n\}, \, \Delta_y(g_i)=\sum_{j=1}^{i-1} \, \ell_{j,i}(y)\, g_j(X) +g_i(y)+\wp(P_i(X,y))$$
where the sum on $j$ is empty for $i=1$ and where $P_i(X,y):= \frac{S(y)^i}{i!}\, X$. 
From $W_i^p-W_i=g_i(X)$, we infer that $\sigma(W_i^p-W_i)=\sigma(g_i(X))=\Delta_y(g_i)$, which implies
$\wp(\sigma(W_i))=\wp( W_i+\sum_{j=1}^{i-1}\, \ell_{j,i}(y) \, W_j+X\, \frac{S(y)^i}{i!} +g_i(y))$. 
Therefore, for all $i$ in $\{1,\cdots, n\}$,
$$\sigma(W_i)=W_i+\sum_{j=1}^{i-1}\, \ell_{j,i}(y) \, W_j+ X\, \frac{S(y)^i}{i!} +Z_i(y) $$
where $Z_i(y)$ is an element of $k$ such that $\wp(Z_i(y))=g_i(y)$.
Using the vector notations of the proposition, we thus obtain the expected formula.
\item To prove the second assertion, we compute $\sigma^p[W]$.
An induction shows that:
$$\sigma^p [W]=(^tL(y))^p [W]+X (\sum_{i=0}^{p-1} \, (^tL(y))^i)\, [\mathcal{R}(y)] \quad \qquad \qquad$$
$$+ y \, (\sum_{i=0}^{p-2} \, (p-i-1)\, (^tL(y))^i )\, [\mathcal{R}(y)]
+ (\sum_{i=0}^{p-1} \, (^tL(y))^i)\, [Z(y)]$$
 We first notice that $(^tL(y))^p$ is equal to the identity matrix $I$,
since $^tL(y)-I$ is nilpotent of size $n \leq p-2$.
Moreover, by Remark 6.2, $^tL(y)= exp (J(y))= \sum_{i=0}^{n-1} \, \frac{J(y) ^i}{i!}$, where $J(y):=S(y) \, ^tJ$.
Accordingly, 
$$
\begin{array}{ll}
\sum_{i=0}^{p-1}\, (^tL(y))^i&=\sum_{i=0}^{p-1} \exp(i\, J(y))\\
&\\
&= I+ \sum_{i=1}^{p-1} \,\sum_{j=0}^{n-1} \frac{(i\,J(y))^j}{j!}=I+\sum_{j=0}^{n-1}\,  \frac{J(y)^j}{j!} \,\sum_{i=1}^{p-1} \,i^j\\
&\\
&=I+(\sum_{i=1}^{p-1}\, i^0) \, I +\sum_{j=1}^{n-1}\,  \frac{J(y)^j}{j!} \,\sum_{i=1}^{p-1} \,i^j\\
&\\
&=\sum_{j=1}^{n-1}\,  \frac{J(y)^j}{j!} \,\sum_{i=1}^{p-1} \,i^j \, \mod \; p
\end{array}
$$
But one easily checks that $\mathcal{N}(j):= \sum_{i=1}^{p-1}\, i^j=0$ mod $p$ for all $j$ in $\{1,\cdots, p-2\}$. Since $n-1 \leq p-3$, we gather that $\sum_{i=0}^{p-1} \, (^tL(y))^i=0$ mod $p$.\\
To conclude, the last sum to compute is $\mathfrak{S}:=\sum_{i=0}^{p-2} \, (p-i-1)\, (^tL(y))^i $. Likewise,
one shows 
$$
\begin{array}{lll}
\mathfrak{S} &=\sum_{i=0}^{p-2}(p-i-1) \exp(iJ(y))&\\
&\\
&=(p-1) \, I + \sum_{i=1}^{p-2} \, (p-i-1) \sum_{j=0}^{n-1} \, \frac{(iJ(y))^j}{j!}&\\
&\\
&= (p-1)\, I- \sum_{i=1}^{p-2} (i+1) \, I - \sum_{j=1}^{n-1} \, \frac{J(y)^j}{j!} \,(\sum_{i=1}^{p-1} \,(i+1)\, i^j) \,& \mod \; p\\
&\\
&= (p-1)\, I -( \mathcal{N}(1)-1) \, I-  \sum_{j=1}^{n-1} \, \frac{J(y)^j}{j!} \, \mathcal{N}(j)- \sum_{j=2}^{n} \, \frac{J(y)^j}{j!} \mathcal{N}(j) \,& \mod \; p
\end{array}$$
Since $\mathcal{N}(j)=0$ when $1 \leq j \leq n \leq p-2$, it follows that $\sum_{i=0}^{p-2} \, (p-i-1)\, (^tL(y))^i=0$ mod $p$.
 As $\sigma^p(X)=X+p\,y=X$ mod $p$, we gather that the order of $\sigma$ divides $p$. Therefore, the group $G[n]$ has exponent $p$. $\square$
\end{enumerate}

\begin{proposition}
Let $(C[p-1],G[p-1])$ be the big action described in the second point of Proposition 6.1, i.e. with 
$n =p-1$. We keep the notations introduced in Proposition 6.1. For all $y$ in $k$, we also define $T(X,y):=\sum_{i=1}^{p-1} \, \frac{(-1)^{i+1}}{i}\, X^i\, y^{p-i}$, i.e. the reduction mod $p$ of $\frac{1}{p} \{(X+\tilde{y})^p-X^p-\tilde{y}^p\} \in W(k)[X]$.
\begin{enumerate}
\item Let $\sigma$ in $G[p-1]$. Let $y$ in $V$ such that $y:=\sigma(X)-X$.
Then,
$$\sigma [W]= ^tL(y) [W]+ X \, [\mathcal{R}(y)]+ [Z(y)]+ [\mathcal{T}(X,y)]$$
where $^tL(y)$ denotes the transpose matrix of the matrix $L(y)$ defined in section 2.4,\\ $[W]:=^t[W_1,\cdots,W_{p-1}]$,  $[\mathcal{R}(y)]:=^t[\frac{S(y)}{1\,!}, \cdots, \frac{S(y)^{p-1}}{(p-1)\, !}]$, $[\mathcal{T}(X,y)]=^t[0,0,\cdots, 0,T(X,y)]$ and $[Z(y)]:=^t[Z_1(y),\cdots,Z_{p-1}(y)]$ where, for all $i$ in $\{1,\cdots,n\}$, $Z_i(y)$ is an element of $k$ satisfying $\wp(Z_i(y))=g_i(y)$.
\item Let $\sigma \in G[p-1]$ as in the first point. Then, if $y \neq 0$, $\sigma$ has order $p^2$. Otherwise, 
the order of $\sigma$ divides $p$. In particular, the group $G[p-1]$ has exponent $p^2$. 
\end{enumerate}
\end{proposition}

\noindent \textbf{Proof:}
\begin{enumerate}
\item The proof of the second point of Proposition 6.1 shows that
$$\forall \, y \in V, \quad \Delta_y(g_{p-1})=\sum_{j=1}^{p-2} \, \ell_{j,p-1}(y)\, W_j +\wp(P_{p-1}(X,y))+g_{p-1}(y)$$
where $$P_{p-1}(X,y):=\frac{S(y)^{p-1}}{(p-1)!} \, X+\sum_{i=1}^{p-1} \, \frac{(-1)^{i+1}}{i}\,
X^i\, y^{p-i}= \frac{S(y)^{p-1}}{(p-1)!} \, X+T(X,y)$$
The same calculation as in the proof of Proposition 6.3 yields:
$$\sigma(W_{p-1})=W_{p-1}+\sum_{j=1}^{p-2} \, \ell_{j,p-1}(y)\, W_j + \frac{S(y)^{p-1}}{(p-1)!} \, X+T(X,y)+Z_{p-1}(y)$$
where $Z_{p-1}(y)$ is an element of $k$ such that $\wp(Z_{p-1}(y))=g_{p-1}(y)$.
The formula of the first point then derives from the vector notation together with the expression
of the others $\sigma(W_i)$, for $1\leq i \leq p-2$, obtained in Proposition 6.3.
\item We now calculate $\sigma^p[W]$.
As in the previous proof, an induction shows that:
$$\sigma^p [W]=(^tL(y))^p [W]+X (\sum_{i=0}^{p-1} \, (^tL(y))^i)\, [\mathcal{R}(y)]
+ y \, (\sum_{i=0}^{p-2} \, (p-i-1)\, (^tL(y))^i )\, [\mathcal{R}(y)]$$
$$ + (\sum_{i=0}^{p-1} \, (^tL(y))^i)\, [Z(y)]+\sum_{i=0}^{p-1} [\mathcal{T}(X+i\,y,y)]$$
Still as in the proof of Proposition 6.3, $^tL(y)^p=I$  and $\sum_{i=0}^{p-1}\, (^tL(y))^i=\sum_{j=1}^{p-2}\,  \frac{J(y)^j}{j!} \,\mathcal{N}(j)$, where $\mathcal{N}(j):=\sum_{i=1}^{p-1} \,i^j=0$ mod $p$, for $j$ in $\{1,\cdots, p-2\}$. 
Besides, as previously seen, $$
\begin{array}{ll}
\sum_{i=0}^{p-2} \, (p-i-1)\, (^tL(y))^i &=- \sum_{j=1}^{p-2} \, \frac{J(y)^j}{j!}\, \{\mathcal{N}(j)+\mathcal{N}(j+1)\}\\
&\\
&=-\frac{J(y)^{p-2}}{(p-2)!} \, \mathcal{N}(p-1)=\frac{J(y)^{p-2}}{(p-2)!}=J(y)^{p-2} \, \mod  \, p
 \end{array}$$ 
Then, $y\,(\sum_{i=0}^{p-2} \, (p-i-1)\, (^tL(y))^i) \, [\mathcal{R}(y)]=y\, (S(y)^tJ)^{p-2} \, [\mathcal{R}(y)]=^t[0,\cdots,0,y\, S(y)^{p-1}]$ mod $p$.
To complete the calculation, one has to compute $\sum_{i=0}^{p-1} T(X+i\,y,y)$.
As $T$ is the reduction mod $p$ of the polynomial $\frac{1}{p} \, \{(X+\tilde{y})^p-X^p-\tilde{y}^p\}$, it follows that 
$$
\begin{array}{lll}
\sum_{i=0}^{p-1} T(X+i\,y,y)&=\frac{1}{p} \sum_{i=0}^{p-1} \{(X+(1+i)\, \tilde{y})^p-(X+i\, \tilde{y})^p- \tilde{y}^p\} \,& \mod \, p \\
&\\
&= \frac{1}{p} \, \{(X+p\,\tilde{y})^p-X^p-p\, \tilde{y}^p\}=-y^p \,& \mod \,p
\end{array}
$$
Therefore, $\sigma[W]^p=[W]+ ^t[0,0,0,\cdots,0,y\,S(y)^{p-1}-y^p]$ mod $p$. 
If $y \in \F_p$, or equivalently $S(y)=0$, then $\sigma[W]^p=[W]+ ^t[0,0,0,\cdots,0,y]$ mod $p$. Otherwise, $S(y) \neq 0$ and, as $S(y)$ lies in $\F_p$, $S(y)^{p-1}=1$ mod $p$, which implies
$\sigma[W]^p=[W]+ ^t[0,0,0,\cdots,0,y-y^p]=[W]+ ^t[0,0,0,\cdots,-S(y)]$ mod $p$.
We gather that if $y=0$, the order of $\sigma$ divides $p$. Otherwise, $\sigma$ has order $p^2$.
This proves the second assertion. $\square$
\end{enumerate}

Contrary to the preceding propositions, the next result describing the center of the group $G[n]$ is common to both cases of Proposition 6.1, i.e. $n<p-1$ and $n=p-1$.
\begin{proposition}
Let $(C[n],G[n])$ be a big action as described in the first or the second point of Proposition 6.1, i.e. with 
$n< p-1$ or $n=p-1$. 
Let $\sigma$ in $G[n]$. Then, $\sigma$ belongs to the center of $G[n]$ if and only if 
\begin{displaymath}
\sigma(X)=X \quad \mbox{and} \quad \forall \, i \in \{1,\cdots,n-1\}, \; \sigma(W_i)=W_i \qquad
\end{displaymath}
 It follows that the center of $G[n]$ is isomorphic to $\Z/p\Z$.
\end{proposition}

\noindent \textbf{Proof:} 
Throughout the proof, we keep the notations introduced in Proposition 6.3 and 6.4.
\begin{enumerate}
\item We first focus on the case $n<p-1$.
Let $\sigma$ in $Z(G[n])$ and $y$ in $V$ such that $y:=\sigma(X)-X$. By definition of $\phi(y)$, for all $g$ in $G[n]_2$, $\phi(y)(g)=\sigma^{-1} \, g \, \sigma=g$, since $\sigma$ lies in the center of $G[n]$. So $\phi(y)=\id$ and $L(y)$ is the identity matrix.
It follows from Remark 6.2 that $S(y)=\ell_{1,2}(y)=0$. Then, $g_i(y)=0$,
for all $1\leq i \leq n$, and $Z_i(y)$, which satisfies $\wp(Z_i(y))=g_i(y)$, lies in $\F_p$. So, by Proposition 6.3, $\sigma[W]=^tL(y) [W]+ X [\mathcal{R}(y)]+ [Z(y)]=[W]+[Z(y)]$.
We now choose some $u$ in $V$ such that $S(u)\neq 0$ and consider $\tau$ in $G[n]$ such that $\tau(X)=X+u$. Then, still by Proposition 6.3,
$\tau[W]=^tL(u)[W]+ X [\mathcal{R}(u)]+ [Z(u)]$. 
One checks that
$\sigma \tau[W]=\tau \sigma[W]$ if and only if $y [\mathcal{R}(u)]+(^tL(u)-I) [Z(y)]=0$.
As $S(u) \neq 0$, it implies that $y=0$ and $Z_i(y)=0$ for all $i$ in $\{1,\cdots,n-1\}$. \\
Conversely, if $y=0$, then $S(y)=0$. As above, it implies $Z_i(y) \in \F_p$, for $1\leq i \leq n$, and $L(y)=I$. It follows that
$\sigma[W]=[W]+[Z(y)]$. Consider $\tau$ in $G[n]$ and $u$ in $V$ such that $u:=\tau(X)-X$. On the one hand,
$\sigma \tau(X)=X+u+y=\tau \sigma(X)$. On the other hand, as seen above, $\sigma \tau[W]=\tau \sigma[W]$ if and only if $y [\mathcal{R}(u)]+(^tL(u)-I) [Z(y)]=0$. If $S(u)=0$, this equality is trivially true. Otherwise, it comes from $Z_i(y)=0$ for all $i$ in $\{1,\cdots,n-1\}$.
Therefore, $\sigma$ lies in the center of $G[n]$.\\
 As a conclusion, $\sigma$ belongs to the center of $G[n]$ if and only if $\sigma(X)=X$,  $\sigma(W_i)=W_i$, with $1 \leq i \leq n-1$,  and $\sigma(W_n)=W_n+Z_n$, with $Z_n$ in $\F_p$. Thus, $Z(G[n]) \simeq \Z/p\Z$.
\item When $n=p-1$, the result is the same but the proof is slightly different. As above, we first notice that, if $\sigma$ lies in $Z(G[n])$, $L(y)=I$,  $S(y)=0$ and then $y$ is in $\F_p$. Write $\tilde{y}^p=\tilde{y}+pR$ with $R$ in $W(k)$. It implies that $g_{p-1}(y)$, which is the reduction mod $p$ of
$\frac{1}{p!}\,(\{\tilde{y}^p-\tilde{y}\}^p-\tilde{y}^{p^2}+\tilde{y}^p)$, is zero. Accordingly, $Z_{p-1}(y)$ also lies in $\F_p$.
Then, by Proposition 6.4, $\sigma[W]=^tL(y) [W]+ X [\mathcal{R}(y)]+ [Z(y)]+[\mathcal{T}(X,y)]=[W]+[Z(y)]+[\mathcal{T}(X,y)]$.
 We now choose some $u$ in $V$ such that $S(u)\neq 0$ and consider $\tau$ in $G[n]$ such that $\tau(X)=X+u$.
 Still by Proposition 6.4, 
$\tau[W]=^tL(u)[W]+ X [\mathcal{R}(u)]+ [Z(u)]+ [\mathcal{T}(X,u)]$. 
One checks that $\sigma \tau[W]=\tau \sigma[W]$ if and only if
\begin{equation} \label{6}
 y \, [\mathcal{R}(u)]+(^tL(u)-I) [Z(y)]+^tL(u)
[\mathcal{T}(X,y)]+[\mathcal{T}(X+y,u)]-[\mathcal{T}(X,u)]-[\mathcal{T}(X+u,y)]=0 
\end{equation}
 As $S(u) \neq 0$, it implies that $y=0$. Then, $T(X,y)=T(X+u,y)=0$ and $T(X+y,u)=T(X,u)$. Thus, one gets: $(^tL(u)-I) [Z(y)]=0$, which implies $Z_i(y)=0$ for all $i$ in $\{1,\cdots,n-1\}$. \\
Conversely, if $y=0$, then $S(y)=0$, $L(y)=I$ and $T(X,y)=0$. It follows that
$\sigma[W]=[W]+[Z(y)]$. Consider $\tau$ in $G[n]$ and $u$ in $V$ such that $u:=\tau(X)-X$. On the one hand, $\sigma \tau(X)=X+u+y=\tau \sigma(X)$. On the other hand, as seen above, $\sigma \tau[W]=\tau \sigma[W]$ if and only if \eqref{6} is satisfied. If $S(u)=0$, this equality is trivially true. Otherwise, it comes from $Z_i(y)=0$ for all $i$ in $\{1,\cdots,n-1\}$ together with $T(X,y)=T(X+u,y)=0$ and $T(X+y,u)=T(X,u)$ obtained because $y=0$.
Therefore, $\sigma$ lies in the center of $G[n]$ and we conclude as in the previous case. $\square$
\end{enumerate}

\begin{corollaire}
Let $p \geq 3$ and $2 \leq n \leq p-1$. We keep the notation of Proposition 6.1. 
\begin{enumerate}
\item The group $G[1]$ is the extraspecial group of order $p^3$ and exponent $p$, namely the unique non abelian group of order $p^3$ and exponent $p$. Moreover, we get the following exact sequence:
$$0 \longrightarrow Z(G[1])\simeq \Z/p\Z \longrightarrow G[1] \longrightarrow (\Z/p\Z)^2\longrightarrow 0$$
\item We also have the following exact sequence:
$$0 \longrightarrow Z(G[n]) \simeq \Z/p\Z \longrightarrow G[n] \longrightarrow G[n-1]\longrightarrow 0$$
\end{enumerate}
\end{corollaire}

\noindent \textbf{Proof:}
\begin{enumerate}
\item The first assertion derives from \cite{LM} (Prop. 8.1).
\item As in Proposition 6.1, we call $K_n:=k(C[n])=k(X,W_1,\cdots, W_n)$ the function field of the curve $C[n]$. Put $k(T):=K_n^{G[n]}$, where $T=Q(X)$, $Q$ being defined as in Proposition 6.1.
 Then, Galois theory, combined with Proposition 6.5, gives the following exact sequence:
$$0 \longrightarrow Gal(K_n/K_{n-1}) \simeq Z(G[n]) \simeq \Z/p\Z \longrightarrow Gal(K_n/k(T)) \longrightarrow Gal(K_{n-1} /k(T)) \longrightarrow 0$$
The claim directly follows.
$\square$
\end{enumerate}

\begin{remarque}
Computation using MAGMA package on finite groups shows that, for $n\geq 2$, the group $G[n]$ is, in general, not uniquely determined by the group extension conditions mentionned in Corollary 6.6.
\end{remarque}

\begin{defi} Following \cite{Ha} and \cite{BT}, we say that a group $G$ is capable if there exists a group $\Gamma$ such that $G \simeq \frac{\Gamma}{Z(\Gamma)}$.
\end{defi}	

\noindent We deduce from Corollary 6.6 the following:
\begin{corollaire}
Let $p \geq 3$ and $2\leq n \leq p-1$. 
\begin{enumerate}
\item The group $G[n-1]$ is capable, with $\Gamma=G[n]$. 
\item In particular, the extraspecial group of order $p^3$ and exponent $p$, with $p >2$, is capable with $\Gamma=G[2]$, a group of order $p^4$. 
\end{enumerate}
\end{corollaire}

\begin{remarque}
Note that \cite{BT} (Example 1.12 p.187) gives another proof of the second statement of Corollary 6.9. Nevertheless, this proof uses some group $\Gamma$ of order $p^5$.
\end{remarque}

\subsubsection{General case.}
Starting from the big actions defined in Proposition 6.1, for which $s_1=1$, we use the base change displayed in \cite{MR} (section 3)
to obtain new ones which still satisfy the conditions of Theorem 5.8 but have arbitrary large $s_1$.

\begin{proposition}
Let $p \geq 3$, $1 \leq n \leq p-1$ and $s_0 \in \N$. Let $S_0(X)$ be an additive separable polynomial of $k[X]$ with degree $p^{s_0}$. Let $(C[n],G[n])$ be the big action defined in Proposition 6.1. Consider the additive polynomial map
$S_0: \p_k^1 \rightarrow C[n]/G[n]_2 \simeq \p_k^1$. 
\begin{enumerate}
\item Let $\tilde{C}[n]:=C[n] \times_{\p_k^1} \p_k^1$
be the curve obtained after the base change defined by $S_0$. 
Then, the cover $\tilde{C}[n] \rightarrow C[n]/G[n]$ is Galois with group 
$\tilde{G}[n] \simeq G[n] \times (\Z/p\Z)^{s_0}$. 
Moreover, the pair $(\tilde{C}[n], \tilde{G}[n])$ is a big action with $\tilde{G}[n]_2\simeq G[n]_2\times \{0\}$ and $Z(\tilde{G}[n])\simeq (\Z /p\Z)^{s_0+1}$. 
\item This big action $(\tilde{C}[n], \tilde{G}[n])$ satisfies the conditions of Theorem 5.8
 with $s_1=s_0+1$.
\end{enumerate}
\end{proposition}

\noindent \textbf{Proof:}
\begin{enumerate}
\item The first assertion derives from \cite{MR} (Prop. 3.1). Another proof consists in replacing $X$ with $S_0(X)$ in the proof of Proposition 6.1, knowing that the calcultation only requires $S_0$ to be additive.
 \item One deduces from Lemma 3.7.2 and Lemma 3.7.7 that $f_i(X) \in \Sigma_{i+1}-\Sigma_i$ implies $f_i(S_0(X)) \in \Sigma_{i+1}-\Sigma_i$. The claim follows. Another proof consists in considering the filtration $(\Lambda_i(G[n]))_{i\geq 0}$, as defined in section 5.1. By Proposition 6.1, this filtration satisfies the first condition of Proposition 5.2. Then, one concludes by checking that, for all $i \geq 0$, $\Lambda_i(\tilde{G}[n])\simeq \Lambda_i(G[n])$.  $\square$
\end{enumerate}

\subsection{A universal family.}
\indent Under the hypotheses of Theorem 5.8, one already knows the form of the functions $f_i$'s, namely their degree $m_i=1+i\, p^{s_1}$ and their belonging to $\Sigma_{i+1}-\Sigma_i$. For given $p$, 
$s_1$ and $n \leq p-1$, this naturally yields an algorithmic method to parametrize the functions $f_i$'s. In this way, we obtain a universal family parametrizing the big actions $(C,G)$ that satisfy Theorem 5.8 with $f_1$ monic and $s_1=1$.
Eventhough it theoretically works for any $p \geq 3$, in what follows, we merely illustrate this method in the special case $p=5$ and $n \leq p-1=4$. 
In this case, we also describe the corresponding space of parameters and, when $n=2$, we give necessary and sufficient conditions on the parameters for two curves of the family to be isomorphic. 
We eventually characterize the subfamily corresponding to the special curves that are studied in section 6.1.1. Throughout this section, the notations concerning big actions are still those fixed in section 3.2.

\begin{proposition}
Fix $p=5$.
Let $(C,G)$ be a big action such that $G_2 \simeq (\Z/p\Z)^n$, with $2 \leq n \leq p-1$. We suppose that $s_1=1$.
We also assume that $(C,G)$ satisfies the conditions of Theorem 5.8. 
Then, there exists 
a coordinate $X$ for the projective line $C/G_2 \simeq \p_1$ and an adapted basis for $A$ as follows:\\

\noindent \underline{For $n=2$:}
$$
\begin{array}{ll}
f_1(X)=X^6+ 2 \, \frac{b_0^{24}+1}{b_0^4} \, X^2\\
\quad \\
f_2(X)= b_0^5 \, X^{11} + 4\,b_0^{25} \, X^7 + 3 \, \frac{4\,b_0^{48} +1}{b_0^3 }\, X^3+b_5\, X\\
\end{array}
$$
Therefore, the parametrization of the functions $f_i$ 's requires two algebraically independent parameters,
namely $b_0$ and $b_5$ in $k$, with $b_0 \neq 0$.\\
\bigskip

\noindent \underline{For $n=3$:}
$$
\begin{array}{ll}
f_1(X)=X^6+ 2 \, \frac{b_0^{24}+1}{b_0^4} \, X^2\\
\quad \\
f_2(X)= b_0^5 \, X^{11} + 4\,b_0^{25} \, X^7 + 3 \, \frac{4\,b_0^{48} +1}{b_0^3 }\, X^3+2\frac{c_7-c_7^5}{b_0^5}\, X\\
\quad \\
f_3(X)= 4\, b_0^{10} \, X^{16} +4\,b_0^{30} \, X^{12} +4\,b_0^{50} \, X^8+ c_7^5\, X^6+
4 \, \frac{b_0^{72}+1}{b_0^2} \, X^4+ 2\,c_7\, \frac{c_7^4\, b_0^{24}+1}{b_0^4} \, X^2+c_9\,X
\end{array}
$$
Thus, the parametrization of the functions $f_i$ 's requires three algebraically independent parameters,
namely $b_0$, $c_7$ and $c_9$ in $k$, with $b_0\neq 0$.\\

\noindent \underline{For $n=4$:}
$$
\begin{array}{ll}
f_1(X)&=X^6+ 2 \, \frac{b_0^{24}+1}{b_0^4} \, X^2\\
\quad \\
f_2(X)&= b_0^5 \, X^{11} + 4\,b_0^{25} \, X^7 + 3 \, \frac{4\,b_0^{48} +1}{b_0^3 }\, X^3+2\frac{c_7-c_7^5}{b_0^5}\, X\\
\quad \\
f_3(X)&= 4\, b_0^{10} \, X^{16} +4\,b_0^{30} \, X^{12} +4\,b_0^{50} \, X^8+ c_7^5\, X^6+
4 \, \frac{b_0^{72}+1}{b_0^2} \, X^4+ 2\,c_7\, \frac{c_7^4\, b_0^{24}+1}{b_0^4} \, X^2+2\frac{d_{11}-d_{11}^5}{b_0^5} X \\
\quad \\
f_4(X)&= 2\,b_0^{15} X^{21} +b_0^{35} \, X^{17} +4\,b_0^{55} \, X^{13}+ d_8^5 \, b_0^5 X^{11} +3 \,b_0^{75}\,X^9 +(4 \,d_8^{25}\, b_0^{25}+4\,b_0^{25}\, c_7^5+b_0^{25}\, c_7^{25}) \, X^7\\
\quad \\
& +d_{11}^5\, X^6+(\frac{b_0^{24}+b_0^{48}}{b_0^3}\, c_7^{25}+ \frac{2+4\, b_0^{24}+4\, b_0^{48}}{b_0^3} \, c_7^5+ \frac{3\,c_7}{b_0^3} +\frac{4\,b_0^{48}+4\,b_{0}^{24}}{b_0^3} \, d_8^{25}+
\frac{b_0^{24}+3+3\,b_0^{48}}{b_0^3}\, d_8^5)
 \, X^3\\
\quad \\
&+2\, \frac{d_{11}\, (d_{11}^4\, b_0^{24}+1)}{b_0^4}\, X^2+d_{13}\, X\\
&\\
\end{array}
$$
with 
$$ b_0^{96}=1 \quad \mbox{and} \quad 2\,t+(3\,b_0^{24} +3) t^5+2\,b_0^{24}\, t^{25}=0 \quad  \mbox{where}
 \quad t:=d_8-c_7$$
 Accordingly, the parametrization of the functions $f_i$ 's requires three algebraically independent
parameters, namely $c_7$, $d_{11}$ and $d_{13}$ in $k$.\\
\end{proposition}

\noindent \textbf{Proof:} We recall that, after an homothety and a translation, one can rigidify the parametrization and fix a coordinate $X$ for the projective line $C/G_2 \simeq \p_1$ such that $f_1$ is a monic polynomial with no monomial of degree one. Furthemore, for $n\geq3$, one also rigidify the functions $f_i$'s by assuming, following Proposition 5.4,
that $\ell_{i,i+1}=\ell_{1,2}$. 
Thus, keeping the writing of exponents in $5$-adic expansion, we write the functions $f_i$'s as follows:
$$
\begin{array}{ll}
f_1(X)&=X^{1+5}+a\, X^2 \\
&\quad \\
f_2(X)&= b_0^5\, X^{1+2.5}+ b_1\, X^{2+5} + b_2\, X^3+ b_3\, X^{1+5} +b_4\,
X^2+b_5\, X\\
&\quad \\
f_3(X)&= c_0\, X^{1+3.5}+c_1\, X^{2+2.5} +c_2\, X^{3+5} +c_3\, X^4+ 
c_4^5\, X^{1+2.5}+ c_5\, X^{2+5}\\
&\\
 &+ c_6\, X^3+ c_7\, X^{1+5} +c_8\,X^2+c_9\, X \\
\quad \\
f_4(X)&= d_0\, X^{1+4.5}+d_1\, X^{2+3.5} +d_2\, X^{3+2.5}+d_3\, X^{4+5}+
b_0^{10}\, d_4\, X^{1+3.5}+d_5\, X^{2+2.5} +d_6\, X^{3+5}\\
&\quad \\
&+d_7\, X^4+ 
b_0^5 \, d_8^5\, X^{1+2.5}+ d_9\, X^{2+5} + d_{10}\, X^3+ d_{11}^5\, X^{1+5} +d_{12}\,
X^2+d_{13}\, X 
\end{array}
$$
with  $b_0 \neq 0$, $c_0 \neq 0$ and  $d_0 \neq 0$.
Note that, for convenience of calculation, some coefficients are directly written as $p$-powers.
Following Proposition 2.13, we first calculate $Ad_{f_1}(Y)=Y^{25}+2\,a^5\, Y^5+Y$. As $V$ is included in $Z(Ad_{f_1})$
and as, in our case, these two vector spaces have the same dimension over $\F_p$, namely $s_1+1=2=2\,s_1$, we gather that $V=Z(Ad_{f_1}).$
We now focus on the relation: 
\begin{equation} \label{7}
\forall \, y \in V, \qquad \Delta_y(f_2)=\ell_{1,2}(y)\, f_1(X) \qquad \mod \; \wp(k[X]) 
\end{equation}
Computations using Maple show that for all $y$ in $V$, $\ell_{1,2}(y)= 2\, b_1\, y +2\, b_0^5 \, y^5$. 
As $V=Z(Ad_{f_1})$, we deduce from Proposition 2.9 that $Ad_{f_1}(X)$ divides the polynomial $(2\,b_1\,X +2\,b_0^5 \, X^5)^5-(2\,b_1\,X +2\,b_0^5 \, X^5)$.
This requires: $b_1= 4\, b_0^{25}$ and $a=2\, \frac{b_0^{24}+1}{b_0^4}$. In addition,
\eqref{7} also yields $b_2 = 3\, \frac{4\, b_0^{48}+1}{b_0^3}$ and $b_3 \in \F_5$. 
Accordingly, by replacing $f_2$ with $f_2-b_3\,f_1$, one can assume that $b_3=0$. It follows that $b_4=0$.
We eventually obtain the expected expression of the functions $f_1$ and $f_2$.\\
\indent Likewise, the case $n=3$ (resp. $n=4$) is solved by studying the relation:
$$
\forall \, y \in V, \qquad \Delta_y(f_3)=\ell_{1,2}(y)\, f_2(X) + \ell_{1,3}(y)\, f_1(X) \qquad \mod \; \wp(k[X])
$$
$$ (\mbox{resp.} \quad 
\forall \, y \in V, \qquad \Delta_y(f_4)=\ell_{1,2}(y)\, f_3(X)+ \ell_{2,4}(y)\, f_2(X) + \ell_{1,4}(y)\, f_1(X) \qquad \mod \; \wp(k[X]) \, ) \,  \square$$

\begin{remarque}
Keeping the notations of Proposition 6.12 with $n=2$, one can show using \cite{LM} (Prop. 3.3) that the two pairs of parameters $(b_0,b_5)$ and $(b'_0, b_5')$  give isomorphic $k$-curves $C$ if and only if 
$$(\frac{b_0'}{b_0})^{24}=1 \quad \mbox{and} \quad b'_5=\pm \, \frac{b'_0}{b_0} \, b_5$$
\end{remarque}

We now emphasize the link with the special family studied in section 6.1. This allows us to characterize the group $G$, at least for $p=5$ and $n <4$. 
\begin{remarque}
\begin{enumerate}
\item Following Proposition 6.11, we apply the linear base change: $X \rightarrow \lambda  \, X$ , with $\lambda \in k^{\times}$,  to the big action defined in Proposition 6.1, for $n \geq 2$. Then, one finds a subfamily of the universal family displayed in Proposition 6.12 if and only if $\lambda \in \F_{25}^{\times}$. For instance, in the case $n=2$, the subfamily  obtained in this case is the one characterized by $b_0^{24}=1$ and $b_5=0$. 
\item We notice that the spaces of parameters of the universal family described in Proposition 6.12 when $n <4$, are Zariski opens of linear affine spaces, which implies that they are irreducible and so connected.
It follows from the preceding point and from Proposition 6.11 (with $s_0=0$) that the group $G$ mentionned in Proposition 6.12 is isomorphic to the one obtained in Proposition 6.1. But, for $n=4$, the space of parameters is no more connected (cf. $b_0^{96}=1$), so the question remains open. 
\end{enumerate}
\end{remarque}

For given $p$ and $2\leq n<p-1$, one does not know the connected components of the space of parameters. Nevertheless, the group structure can be approached via the proposition below that generalizes Proposition 6.5 and Corollary 6.6.
\begin{proposition}
Let $p \geq 3$ and $2 \leq n \leq p-1$. Let $(C,G)$ be a big action which satisfy the conditions of Theorem 5.8 with $s_1=1$. 
\begin{enumerate}
\item The center of $G$, $Z(G)$, is included in its derived subgroup $D(G)$. It follows that $Z(G)$ is cyclic of order $p$. 
\item Moreover, for all $i$ in $\{1, \cdots, n\}$, the quotient group $G/ \Lambda_i(G)$ is capable.
\end{enumerate}
\end{proposition}

\noindent \textbf{Proof:}
\begin{enumerate}
\item  As $G$ satisfies the conditions of Theorem 5.8, $\Lambda_{n-1}(G)$ is an index $p$ subgroup of $G_2=D(G)$.
As $\Lambda_{n-1}(G)=(G_2)^{\nu_1}$ (cf. Theorem 5.8), the quotient curve $C/\Lambda_{n-1}(G)$ 
is the $p$-cyclic cover of the affine line parametrized by $W_1^p-W_1=f_1(X)$. Since $v=s_1+1=2$, it follows from \cite{LM} that  the group $G/\Lambda_{n-1}(G)$ is the extraspecial group of order $p^3$ and exponent $p$. In particular, its center is a $p$-cyclic group generated by $\tau$
such that $\tau(X)=X$ and $\tau(W_1)=W_1+1$. 
Now, take $\sigma \in Z(G)$. Then, $\sigma$ induces $\tilde{\sigma} \in Z(G/\Lambda_{n-1}(G))$. So, $\sigma(X)=X$. As $k(X)=L^{D(G)}$, it implies that $Z(G)$ is included in $D(G)$.
Besides, by Theorem 5.8, $\Lambda_1(G)= Z(G)\cap D(G)=Z(G)\simeq \Z/p\Z$, which proves that $Z(G)$ is $p$-cyclic.
\item  Theorem 5.8 implies that the cover $C \rightarrow C/G_2$ is parametrized by $n$ Artin-Schreier equations: $W_j^p-W_j=f_j(X) \in \Sigma_{j+1}-\Sigma_j$, with $1 \leq j \leq n$. Take $i$ in $\{0, \cdots, n-1\}$. Then, the curve $C/\Lambda_i(G)$ is parametrized by the $n-i$ first equations: $W_j^p-W_j=f_j(X) \in \Sigma_{j+1}-\Sigma_j$, with $1 \leq j \leq n-i$. It follows that the pair $(C/\Lambda_i(G), G/ \Lambda_i(G))$ is a big action (cf. \cite{MR} Lemma 2.4) which still satisfies Theorem 5.8 with $s_1=1$. 
We deduce from the first point that $\Lambda_1(G/\Lambda_i(G))=Z(G/\Lambda_i(G))$. As 
$\frac{G/\Lambda_i(G)}{\Lambda_1(G/\Lambda_i(G))} \simeq G/\Lambda_{i+1}(G)$, we get the exact sequence: 
$$0 \longrightarrow Z(G/\Lambda_{i}(G)) \longrightarrow G/\Lambda_{i}(G) \longrightarrow G/\Lambda_{i+1}(G)\longrightarrow 0$$ The claim follows. $\square$

\end{enumerate}

\noindent We conclude with the following 

\medskip

\noindent \textbf{Problems:} 
\begin{enumerate}
\item For any $p$, find equations for the universal family (at least for $s_1=1$) as we obtained for the special family.
\item Compare the universal family corresponding to a given $s_1$ with the one obtained after a base change by a generic and additive polynomial map, applied to the universal family with $s_1=1$.
\end{enumerate}

\noindent A last interesting question is raised by the following
\begin{remarque}
Proposition 6.12 seems to suggest that any $p$-cyclic \'etale cover of the affine line given by $$W_1^p-W_1=f_1(X):=X\,S(X) \quad  \mbox{with} \quad S \in k\{F]$$  could be embedded in a big action $(C,G)$ where $C$ is parametrized by $n$ Artin-Schreier equations: 
$$W_i^p-W_i=f_i(X) \in \Sigma_{i+1}-\Sigma_i \quad \mbox{with} \quad 1 \leq i \leq n <p-1$$ and that without any restriction on the coefficients of $f_1$. Nevertheless, it is no more true for $n=p-1$ unless the coefficients of $S(X)$ satisfy a specific algebraic condition to be determined (see
e.g. $b_0^{96}=1$ in Proposition 6.12).
\end{remarque}

\bigskip

\noindent {\bf \Large{Acknowledgements}.}\\
I wish to express my gratitude to my supervisor, M. Matignon, for his valuable advice, his precious help and his permanent support all along the elaboration of this paper.

\bigskip

\begin{flushleft}

Magali ROCHER\\
Institut de Math\'ematiques de Bordeaux,
Universit\'e de Bordeaux I,
351 cours de la Lib\'eration, 
33405 Talence Cedex, France. \\
e-mail : {\tt Magali.Rocher@math.u-bordeaux1.fr}

\end{flushleft}

\end{document}